\documentclass[10pt,a4paper]{amsart}
\usepackage[a4paper,margin=25mm]{geometry}
\usepackage{graphicx} 

\usepackage{verbatim}
\usepackage{amsmath}
\usepackage{amssymb}
\usepackage{cite}
\usepackage{bm}
\usepackage{comment}
\usepackage{hyperref}
\usepackage{overpic}
\usepackage[T1]{fontenc}
\usepackage{lmodern}
\usepackage{microtype}
\theoremstyle{plain}
\newtheorem{thm}{Theorem}[section] 
\theoremstyle{definition}
\newtheorem{defn}[thm]{Definition} 
\theoremstyle{definition}
\newtheorem{defnprop}[thm]{Definition-Proposition}
\theoremstyle{plain}
\newtheorem{cor}[thm]{Corollary} 
\theoremstyle{plain}
\newtheorem{prop}[thm]{Proposition}
\theoremstyle{plain}
\newtheorem{lem}[thm]{Lemma}
\theoremstyle{plain}
\newtheorem{conj}[thm]{Conjecture}
\theoremstyle{remark}
\newtheorem{rem}[thm]{Remark}
\theoremstyle{remark}

\theoremstyle{definition}

\theoremstyle{definition}

\title{Skein and cluster algebras of punctured surfaces}
\author{Enhan Li}
\address{Department of Mathematics, University of Science and Technology of China}
\email{lablehustc0719@mail.ustc.edu.cn}
\date{}
\begin{document}

	\begin{abstract}
    We prove the full Fock--Goncharov conjecture for the $\mathcal{A}_{SL_2,\Sigma_{g,p}}$-cluster variety associated with the moduli space of decorated twisted $SL_2$-local systems on triangulable surfaces $\Sigma_{g,p}$ with at least two punctures. Equivalently, we show that the tagged skein algebra $Sk^{ta}(\Sigma)$, identified with the middle cluster algebra $\mathrm{mid}(\mathcal{A})$, coincides with the upper cluster algebra $U(\Sigma)$. We also introduce a localized cluster variety $\mathring{\mathcal{A}}$ whose positive real locus is the decorated Teichm\"uller space $\mathcal{T}^d(\Sigma)$. For $p\geq2$, its ring of global sections is the normalization of the classical Roger--Yang skein algebra $Sk^{RY}_{q\to1}(\Sigma)$ and equals that algebra except for the four-punctured sphere. In the equality cases, the quantum Roger--Yang skein algebra $Sk^{RY}_q(\Sigma)$ provides a quantization of this global-section ring. We prove normality of $Sk^{ta}(\Sigma)$ for $p\geq2$, complementing its known Gorenstein property, and normality of $Sk^{RY}_{q\to1}(\Sigma)$ for every triangulable $\Sigma_{g,p}$ other than $\Sigma_{0,4}$. We also prove finite generation of $U(\Sigma_{g,p})$ for every triangulable $\Sigma_{g,p}$, including the once-punctured case.
	\end{abstract}

\maketitle

\tableofcontents

		\section{Introduction}

        Let $\Sigma=\Sigma_{g,p}$ be a surface of negative Euler characteristic obtained by removing $p>0$ punctures from a compact orientable surface $\Sigma_g$ of genus $g$. Let $\mathcal A$ denote the associated cluster variety, obtained by gluing seed tori as in Section 2. In the present paper, we focus on the map $\mathcal{X}_{PGL_2,\Sigma}(\mathbb{Z}^t)\to \Gamma(\mathcal A,\mathcal O_{\mathcal A})$ sending an integral $\mathcal{X}$-lamination to a function (in particular, the trace of the monodromy representation for a simple loop) on the moduli $\mathcal{A}_{SL_2,\Sigma_{g,p}}$ of decorated twisted $SL_2$-local systems, as constructed by Fock and Goncharov\cite{fock06moduli}. The vector space spanned by the functions has the structure of an algebra, which turns out to be the tagged skein algebra $Sk^{ta}(\Sigma)$. It is conjectured that the ring of global sections of the cluster variety $\mathcal A$ is spanned by these functions associated to integral $\mathcal{X}$-laminations, so the tagged bracelet elements of $\Sigma_{g,p}$ parametrize a canonical basis on $\mathcal A$.\par

        The study of canonical bases on moduli of local systems (or more specifically, cluster varieties) leads to the interaction of skein and cluster algebras\cite{muller16skein}\cite{ishibashi23au}. As for the moduli space $\mathcal{A}_{SL_2,\Sigma}$ of decorated twisted local systems on $\Sigma$, there are mainly two kinds of skein algebras we are concerned about, namely the tagged skein algebra $Sk^{ta}(\Sigma)$ and the classical (resp. quantum) Roger--Yang skein algebra $Sk^{RY}_{q\to1}(\Sigma)$ (resp. $Sk^{RY}_q(\Sigma)$). Our goal is to compare the tagged skein algebra $Sk^{ta}(\Sigma)$ with the upper cluster algebra $U(\Sigma)$\cite{berenstein05upper}, which is exactly the ring of global sections of the $\mathcal{A}_{SL_2,\Sigma}$-cluster variety. However, instead of directly comparing the two algebras, it turns out to be more natural to compare their geometry. In other words, we establish a relation between the cluster variety $\mathcal{A}$ and the spectrum $\mathrm{Spec}(Sk^{ta}(\Sigma))$ of the tagged skein algebra. We state our main result as follows. It follows from such a geometric description that for a triangulable surface $\Sigma=\Sigma_{g,p}$ with $p\geq 2$, the tagged skein algebra coincides with the upper cluster algebra. 

        \begin{thm}[Theorem \ref{mainthmproof}, Theorem \ref{oncepuncturedconj}]
        (1) For a triangulable surface $\Sigma=\Sigma_{g,p}$ with $p\geq 2$, there is an open immersion $\mathcal{A}\to\mathrm{Spec}(Sk^{ta}(\Sigma))$ from the cluster variety to the spectrum of the tagged skein algebra, which is a morphism between integral Noetherian schemes. The complement of the image of $\mathcal{A}$ has codimension at least 2 in $\mathrm{Spec}(Sk^{ta}(\Sigma))$.\\
        (2) For a triangulable surface $\Sigma=\Sigma_{g,1}$, there is an open immersion $\tilde{\mathcal{A}}\to\mathrm{Spec}(Sk^{ta}(\Sigma))$ from the enlarged cluster variety to the spectrum of the tagged skein algebra, which is a morphism between integral Noetherian schemes. The complement of the image of $\tilde{\mathcal{A}}$ has codimension at least 2 in $\mathrm{Spec}(Sk^{ta}(\Sigma))$.
        \end{thm}

        \begin{cor}[Theorem \ref{mainthmproof}]
         For a triangulable surface $\Sigma=\Sigma_{g,p}$ with $p\geq 2$, the tagged skein algebra $Sk^{ta}(\Sigma)$ coincides with the upper cluster algebra $U(\Sigma)$. As a result, both the tagged bracelets and the theta functions parametrize a basis of the ring of global sections $\Gamma(\mathcal{A},\mathcal{O}_\mathcal{A})$ of the $\mathcal{A}$-cluster variety.
        \end{cor}

        When a spin structure is fixed on $\Sigma$, the set of positive real points $\mathcal{A}_{SL_2,\Sigma}(\mathbb{R}_{>0})$ on the cluster variety recovers Penner's decorated Teichm\"uller space $\mathcal{T}^d(\Sigma)$\cite{fock06moduli}. To extend the skein relations to curves pointing at punctures, Roger and Yang define a lambda function on $\mathcal{T}^d(\Sigma)$ for each simple curve on $\Sigma$\cite{roger14skein}, which is (after a suitable choice of orientation at punctures) exactly the function defined by Fock and Goncharov under the identification $\mathcal{A}_{SL_2,\Sigma}(\mathbb{R}_{>0})=\mathcal{T}^d(\Sigma)$. Roger and Yang define a map $Sk^{RY}_{q\to1}(\Sigma)\to C^{\infty}(\mathcal{T}^d(\Sigma))$ from the classical Roger--Yang skein algebra to the ring of smooth functions on the decorated Teichm\"uller space. Since this map turns out to be injective\cite{moon24consequences}, it is natural to seek a suitable algebraic structure on $\mathcal{T}^d(\Sigma)$ which has $Sk^{RY}_{q\to1}(\Sigma)$ as its ring of global sections\cite{roger14skein}. We give an algebraic structure on $\mathcal{T}^d(\Sigma)$ in terms of the localized cluster variety $\mathring{\mathcal{A}}$, which is the largest open subscheme of $\mathcal{A}$ making the potentials (or horocycle length functions) $v_i^{\pm1}$ regular. For $p\geq2$ and $\Sigma\ne\Sigma_{0,4}$, the quantum Roger--Yang skein algebra then provides a quantization of the coordinate algebra of this localized cluster variety, whose positive real locus is the decorated Teichm\"uller space.
        \begin{thm}[Corollary \ref{localized}]
        For a triangulable surface $\Sigma=\Sigma_{g,p}$ with $p\geq2$, the positive real locus of the localized cluster variety is
        $\mathring{\mathcal A}_{\mathbb R}(\mathbb R_{>0})=\mathcal T^d(\Sigma)$.
        If $\Sigma\ne\Sigma_{0,4}$, its ring of global sections is the classical Roger--Yang skein algebra over $\mathbb R$. For $\Sigma=\Sigma_{0,4}$, the corresponding complex ring of global sections is the normalization of $Sk^{RY}_{q\to1}(\Sigma)$ and is strictly larger.
        \end{thm}

        We remark that analogous results to our main theorem, in the bordered cases and coefficient conventions recalled in Section 2, follow from the cluster theoretic equality $A(\Sigma)=U(\Sigma)$ between the cluster and upper cluster algebras, which unfortunately fails when $\Sigma=\Sigma_{g,p}$ with $g\geq 1$. This compels us to compare the tagged skein algebra $Sk^{ta}(\Sigma)$ and the upper cluster algebra $U(\Sigma)$ without the aid of the (ordinary) cluster algebra $A(\Sigma)$.\par

        As a result, our work combines both cluster and skein theoretic methods. On the cluster side, our geometric description is based on results on acyclic cluster algebras by Muller\cite{muller13acyclic,muller14au} and on once-marked surfaces by Canakci, Lee and Schiffler\cite{canakci15cluster}. To pass from the geometric to the algebraic side, the S2 property of the tagged skein algebra $Sk^{ta}(\Sigma)$ is due to Gross--Hacking--Keel--Kontsevich's construction of a theta basis\cite{gross18canonical}, Mandel--Qin's comparison between bracelet and theta bases\cite{mandel23bracelet}, and Goncharov--Shen's work stating that the Donaldson--Thomas transformation for $\mathcal{X}_{PGL_2,\Sigma_{g,p}}$ is cluster when $p\geq 2$\cite{goncharov18dt}. On the skein side, finite generation was established for the Roger--Yang arc algebra by Bobb, Kennedy, Wong and Peifer\cite{bobb16finite}; the corresponding results for tagged and generalized skein algebras\cite{karuo25cluster,bloomquist23degenerations} are foundational to our methods, and the idea of relating the dimension of skein algebras to that of the space of measured laminations has appeared in the work of Bloomquist--Karuo--L\^e\cite{bloomquist23degenerations}, based on the expression of Dehn--Thurston coordinates in terms of intersection numbers due to Luo and Stong\cite{luo04dt}.

        The Roger--Yang comparison also requires an independent normality argument. We use an invariant-theoretic model related to the decorated character algebras of Muller and Samuelson\cite{muller13decorated}, with curve evaluations and tangent-lift conventions from Costantino and L\^e\cite{costantino22stated}. For the four-punctured sphere, the presentation of Azad, Chen, Dreyer, Horowitz and Moon\cite{azad21presentations} gives an explicit obstruction to normality. This accounts for the normalization appearing in that exceptional case.

		\subsection{Fock--Goncharov duality conjecture}
		
		As an algebro-geometric approach towards higher Teichm\"uller theory, Fock and Goncharov \cite{fock06moduli} introduced the moduli space of decorated twisted local systems $\mathcal{A}_{G,\Sigma}$ and the moduli space of framed local systems $\mathcal{X}_{G,\Sigma}$ for a surface $\Sigma$ with negative Euler characteristic and a split reductive algebraic group $G$ (over $\mathbb{Q}$), known as the Fock--Goncharov $\mathcal{A}$ and $\mathcal{X}$ moduli spaces. When $\Sigma$ is a compact orientable surface minus $p$ punctures, the positive real points on the moduli space $\mathcal{A}_{SL_2,\Sigma}$ (resp. $\mathcal{X}_{PGL_2,\Sigma}$) recover the classical decorated (resp. enhanced) Teichm\"uller spaces. Remarkably, for a connected and simply connected split semi-simple (resp. adjoint split semi-simple) algebraic group, the moduli space $\mathcal{A}_{G,\Sigma}$ (resp. $\mathcal{X}_{G,\Sigma}$) is equipped with a cluster (resp. dual cluster) structure\cite{le17cluster}\cite{goncharov25quantum}, which is a set of toric charts with positive transition functions that are encoded completely in an initial combinatorial object called a quiver.\par
        
        For a connected and simply connected split semi-simple algebraic group $G$ with Langlands dual group $^{L}G$, Fock and Goncharov conjectured  the spaces $\mathcal{X}_{^{L}G,\Sigma}$ and $\mathcal{A}_{G,\Sigma}$ to be dual to each other, in the sense that the integral tropical points on $\mathcal{A}_{G,\Sigma}$ (resp. $\mathcal{X}_{^{L}G,\Sigma}$) parametrize a linear basis of global sections on the dual variety $\mathcal{X}_{^{L}G,\Sigma}$ (resp. $\mathcal{A}_{G,\Sigma}$). Such a parametrization is expected to be equivariant under the action of the mapping class group and to satisfy certain positivity conditions. This is known as the Fock--Goncharov duality conjecture for the moduli space $\mathcal{X}_{^{L}G,\Sigma}$ (resp. $\mathcal{A}_{G,\Sigma}$). \par 

        For $G=SL_2$ and $\Sigma=\Sigma_{g,p}$, an explicit parametrization in the duality conjecture has been given in \cite[Section 12]{fock06moduli}. The set of integer tropical points $\mathcal{X}(\mathbb{Z}^t)$ is interpreted as the set of $\mathcal{X}$-laminations on $\Sigma$, or collections of non-intersecting nontrivial non-peripheral curves with integral weights, together with an orientation at punctures of positive weight. One may choose a spin structure on $\Sigma$, equivalently a splitting of the central $\mathbb Z/2\mathbb Z$-extension $\pi_1(T\Sigma-\Sigma)/\langle z^2\rangle$ of $\pi_1(\Sigma)$, where $z$ is the tangent-fiber generator, to identify twisted local systems with ordinary local systems. Curve evaluations use tangent-framed lifts. Then, a weight-$k$ loop is associated with the trace of the $k$th power of the monodromy along the lifting $\tilde{l}$ to the punctured tangent bundle $T\Sigma-\Sigma$, and a simple arc (with the orientation coherent with that of punctures) is associated with the edge function $\Delta(l)$ pairing the two monodromy-invariant vectors at two punctures through parallel transport along $\tilde{l}$. The construction is extended to general $\mathcal{X}$-laminations by the Weyl group action and the additivity of compatible $\mathcal{X}$-laminations. The parametrization map $\mathcal{X}_{PGL_2,\Sigma}(\mathbb{Z}^t)\to \Gamma(\mathcal A,\mathcal O_{\mathcal A})$ induces a multiplicative pairing between the spaces $\mathcal{A}_{SL_2,\Sigma}$ and $\mathcal{X}_{PGL_2,\Sigma}$, which is intimately related to the hyperbolic length functions. After tropicalization, the canonical pairing  recovers the intersection pairing between $\mathcal{A}$- and $\mathcal{X}$-laminations\cite{fock06moduli}\cite{shen23intersections}. \par

        Fock and Goncharov's construction of the parametrization map $\mathcal{X}_{PGL_2,\Sigma}(\mathbb{Z}^t)\to \Gamma(\mathcal A,\mathcal O_{\mathcal A})$ preserves the skein relations, so the association $\mathcal{X}_{PGL_2,\Sigma}(\mathbb{Z}^t)\to \Gamma(\mathcal A,\mathcal O_{\mathcal A})$ induces a homomorphism from the tagged skein algebra $Sk^{ta}(\Sigma)$ to the ring of global sections $\Gamma(\mathcal A,\mathcal O_{\mathcal A})$ of the cluster variety. As a result of preceding works, the homomorphism $Sk^{ta}(\Sigma)\to \Gamma(\mathcal A,\mathcal O_{\mathcal A})$ turns out to be injective and, for $p\geq2$, the duality conjecture is reduced to the surjectivity of the map (Proposition \ref{conjs}). 

        Fock--Goncharov's $\mathcal{X}$- and $\mathcal{A}$-spaces, as well as the duality conjecture, can be generalized to the setting of cluster ensembles\cite{fock09cluster} and hence to arbitrary triangulable surfaces with boundary marked points. The generalized setting enriches the theory of Fomin--Zelevinsky cluster algebras\cite{fomin01cluster} from a geometric point of view; in particular, the ring of global sections $\Gamma(\mathcal{A},\mathcal{O}_{\mathcal{A}})$ is exactly the upper cluster algebra $U(\Sigma)$. This geometric point of view serves as the foundation of the approach of the present paper to compare skein and upper cluster algebras in terms of their geometry.\par

        On the skein side, bracelet bases and their positivity were developed by Musiker, Schiffler and Williams\cite{musiker13bases} and Thurston\cite{thurston14positive}. Gross, Hacking, Keel and Kontsevich\cite{gross18canonical} constructed theta functions, and Mandel and Qin\cite{mandel23bracelet} identified the tagged bracelet and theta bases in the cases used here. However, for closed punctured surfaces $\Sigma_{g,p}$ in general, these results do not by themselves establish that the middle cluster algebra equals $U(\Sigma)$. For $p\geq2$, $Sk^{ta}(\Sigma)$ is identified with the image of $\mathrm{mid}(\mathcal A)$ in $U(\Sigma)$, while the remaining surjectivity is the main comparison proved in this paper. Davison and Mandel\cite{davison21quantum} developed quantum theta bases and proved their strong positivity. The quantum statements require their specified coefficient and compatibility conventions.

        We give a review of preceding works on skein and cluster algebras and their compatibility in Section 2. Roughly speaking, for a triangulable surface $\Sigma$ (possibly with boundary and marked points), there are canonical inclusions $A(\Sigma)\subseteq Sk^{ta}(\Sigma)\subseteq U(\Sigma)$ where $A(\Sigma)$ is the (ordinary) cluster algebra, and the second inclusion recovers Fock and Goncharov's construction when $\Sigma=\Sigma_{g,p}$. For the bordered cases and coefficient conventions recalled in Section 2, the known equality $A(\Sigma)=U(\Sigma)$ makes the second inclusion an equality and gives the corresponding duality statement. However, for closed punctured surfaces of positive genus, $A(\Sigma)\ne U(\Sigma)$: Ladkani\cite{ladkani13once} established this for the once-punctured case, together with non-finite generation of $A(\Sigma_{g,1})$, and Moon and Wong\cite{moon24consequences} proved the comparison more generally. The cluster structure is also more complicated combinatorially, for instance because local acyclicity fails. We give a proof of the equality $Sk^{ta}(\Sigma)=U(\Sigma)$ for $p\geq 2$ in Section 5, with the aid of both skein and cluster theoretical techniques. For $\Sigma_{g,1}$, the two tag components lead to an enlarged cluster variety. Zhou\cite{zhou20cluster} studied an enlarged atlas for the once-punctured torus, and Section 7 compares the enlarged atlas with the normalization of the tagged skein algebra. Finite generation of $U(\Sigma_{1,1})$ was already known from the explicit presentation of Matherne and Muller\cite{matherne15upper}; Section 7 gives a uniform normalization description for all $g\geq1$.\par

        Our main theorem is exactly the geometric analog of the equality $Sk^{ta}(\Sigma)=U(\Sigma)$. We recall the construction of $\mathcal{A}$ as a union of cluster charts \cite{fock09cluster}\cite{gross15birational}, each included as an open subscheme into $\mathrm{Spec}(Sk^{ta}(\Sigma))$ due to the fact that each tagged arc or loop is a Laurent polynomial in the cluster variables of each seed\cite{roger14skein}\cite{mandel23bracelet}. Now, the inclusion $Sk^{ta}(\Sigma)\to U(\Sigma)$ is exactly the induced map between rings of global sections of the open immersion $\mathcal{A}\to \mathrm{Spec}(Sk^{ta}(\Sigma))$, which will be shown to be an isomorphism up to codimension 2 by Theorem \ref{mainthmproof}.\par

\subsection{Decorated Teichm\"uller space and Roger--Yang's program}
		
	Fock--Goncharov showed in \cite{fock06moduli} that the space of positive real points $\mathcal{A}_{SL_2,\Sigma}(\mathbb{R}_{>0})$ on the moduli space of decorated twisted $SL_2$-local systems is precisely Penner's decorated Teichm\"uller space $\mathcal{T}^d(\Sigma)$ \cite{penner87decorated}, which is the $\mathbb{R}_{>0}^p$-principal bundle over the cusped Teichm\"uller space with each point on the fiber indicating a collection of horocycles at punctures. \par

    Roger and Yang introduced a quantized skein algebra $Sk^{RY}_q$ \cite{roger14skein}, whose specialization $Sk^{RY}_{q\to1}$ when $q\to 1$ admits a homomorphism $\Phi:Sk^{RY}_{q\to1}\to C^{\infty}(\mathcal{T}^d(\Sigma))$ which is consistent with the comparison between $\mathcal{A}_{SL_2,\Sigma}$ and $\mathcal{T}^d(\Sigma)$. More specifically, a curve system $\alpha$ is mapped to the lambda function $\lambda(\alpha)$ generalizing the lambda length of geodesic arcs between punctures and the lifted $SL_2$-trace function of simple loops. The potential $v_i$ is mapped to the length function of the horocycle. Positive potentials on moduli spaces of decorated local systems also play a central role in the canonical-basis and mirror-symmetry framework of Goncharov and Shen\cite{goncharov15mirror}. The map $\Phi$ is an analog of Fock and Goncharov's map $\mathcal{X}_{PGL_2,\Sigma}(\mathbb{Z}^t)\to \Gamma(\mathcal{A}_{SL_2,\Sigma},\allowbreak\mathcal{O}_{\mathcal{A}_{SL_2,\Sigma}})$, which is related to the tagged skein algebra $Sk^{ta}(\Sigma)$. Tagging an endpoint is represented by multiplying its underlying arc by the corresponding potential. The distinction between the two classical algebras is that $Sk^{RY}_{q\to1}(\Sigma)=Sk^{ta}(\Sigma)[v_1^{\pm1},\ldots,v_p^{\pm1}]$ adjoins all potentials and their inverses; these functions are units on its spectrum. In Section 3, we shall give a clear geometric description in terms of the spectra of the two algebras.
		\begin{prop}[Proposition \ref{RYtaggedbirational}]
             Let $\mathcal S$ be the largest open subscheme of $\operatorname{Spec}(Sk^{ta}(\Sigma))$ on which all $v_i^{\pm1}$ are regular. The natural morphism
             $\phi:\operatorname{Spec}(Sk^{RY}_{q\to1}(\Sigma))\to\operatorname{Spec}(Sk^{ta}(\Sigma))$
             induces an isomorphism $\phi^{-1}(\mathcal S)\cong\mathcal S$. Thus $\mathcal S$ is canonically an open subscheme of both spectra.
        \end{prop}

		Via the homomorphism $\Phi:Sk^{RY}_{q\to1}\to {C}^{\infty}(\mathcal{T}^d(\Sigma))$, Roger and Yang attempted to interpret $Sk^{RY}_q$ as a deformation quantization of the decorated Teichm\"uller space $\mathcal{T}^d(\Sigma)$. The work of Moon and Wong\cite{moon21skein}\cite{moon24consequences} and Bloomquist--Karuo--L\^e\cite{bloomquist23degenerations} has completed the first step, which is the injectivity of $\Phi$. Further, $\mathcal{T}^d(\Sigma)$ is expected to be equipped with a suitable algebraic structure in terms of $Sk^{RY}_{q\to1}(\Sigma)$, just as the $SL_2$-character variety $\mathcal{X}(\Sigma)$ for $\Sigma=\Sigma_g$ has the structure of an affine variety whose ring of global sections is generated by trace functions; the classical skein--character correspondence was established by Bullock\cite{bullock97characters} and Przytycki--Sikora\cite{przytycki00characters}. Inspired by Shen--Sun--Weng\cite{shen23SLn}, we shall construct in Section 6 the localized cluster variety $\mathring{\mathcal{A}}_\mathbb{R}$ over $\mathbb{R}$ whose complexification $\mathring{\mathcal{A}}$ is an open subscheme of the usual cluster variety  $\mathcal{A}$. It will follow immediately from our construction that $$\mathring{\mathcal A}_{\mathbb R}(\mathbb R_{>0})=\mathcal T^d(\Sigma).$$ Also, the main theorem in the present paper will imply that for $\Sigma_{g,p}$ with $p\geq 2$ and $(g,p)\ne(0,4)$,  $$Sk^{RY}_{q\to1}=\Gamma(\mathring{\mathcal{A}},\mathcal{O}_{\mathring{\mathcal{A}}})$$
		as conjectured in \cite{shen23SLn}.
		Therefore, combined with the injectivity of $\Phi$ \cite[Theorem A]{moon24consequences}, we obtain this realization of Roger--Yang's program outside the four-punctured sphere, where normalization is necessary.\par 

        On the $\mathcal X$-side, Bonahon and Wong\cite{bonahon11quantum} constructed the quantum trace relating Kauffman bracket skein algebras to quantum Teichm\"uller theory. Allegretti and Kim\cite{allegretti17duality} used it to construct a quantum duality map from tropical $\mathcal A$-points to quantum functions on $\mathcal X$. Ishibashi and Karuo\cite{ishibashi24duality} extended this construction to general marked surfaces and studied compatibility with the dual map when there are no interior marked points. The quantum trace construction concerns the duality direction opposite to the principal comparison in the present paper.

        The coefficient-free exchange matrix for a closed punctured surface is not of full rank\cite{fomin08cluster1}, so it does not admit a compatible matrix in the standard quantum-cluster construction\cite{berenstein05quantum}. Inspired by \cite{shen23SLn}, one may instead consider the subalgebra $Sk^{ta}_q(\Sigma)$ of $Sk^{RY}_q(\Sigma)$ generated by loops and arcs weighted by the puncture variables according to their tags. Identifying its specialization with $Sk^{ta}(\Sigma)$ and establishing a compatible quantum bracelet basis require separate arguments. Karuo\cite{karuo24positivity} studies bases and a positivity conjecture for Roger--Yang skein algebras, explicitly distinguishing that setting from quantum tagged bracelets with cluster coefficients. Our theorem identifies the classical algebra $Sk^{ta}(\Sigma)$ with $U(\Sigma)$; it does not by itself establish this additional quantum tagged-skein comparison.

        \subsection{(Lack of) local acyclicity and estimation of Krull dimension}
        Our proof of the main theorem follows the cutting and freezing principle, which is a well-known technique that has been used in proving properties of locally acyclic algebras\cite{muller13acyclic,muller14au,muller16skein}. We underline the importance of viewing the inclusion $Sk^{ta}(\Sigma)\subseteq U(\Sigma)$ geometrically as the pull-back of the open immersion $\mathcal{A}\to \mathrm{Spec}(Sk^{ta}(\Sigma))$, since we shall frequently localize at cluster variables and compare the localized open immersion with that given by the frozen cluster variety. The operation of freezing is much more natural from the geometric perspective in terms of the cluster variety, which is just given by a subset of cluster charts.\par
        Since the algebraic and combinatorial structure of the cluster variety for $\Sigma=\Sigma_{g,p}$ can be much more subtle than that for most marked surfaces (for instance, covering pairs do not exist anymore\cite{muller13acyclic}), it is not expected that we may cover (the spectrum of) any cluster algebra (e.g. $A(\Sigma)$, $U(\Sigma)$, or $Sk^{ta}(\Sigma)$ which is in the middle) by (the spectra of) acyclic cluster algebras (e.g. those obtained by freezing some cluster variables) completely. However, we shall see that away from a closed subscheme cut out by some ideal $I_0$, the rest of $\mathrm{Spec}(Sk^{ta}(\Sigma))$ can be covered by (spectra of) frozen cluster algebras for which the cluster and upper cluster algebras coincide. This reduces the comparison between skein and upper cluster algebras to the characterization of the ideal $I_0$, especially its height.\par The Louise property used for punctured cut components is independent of the frozen rows, as in Muller and Speyer\cite{mullerspeyer16louise}.

        Matherne and Muller\cite{matherne15upper} established a general criterion recognizing an intermediate Noetherian ring $A\subseteq R\subseteq U$ as the upper cluster algebra using $S_2$ and a codimension-two condition. The surface-specific input here is to obtain the required dimension bounds by filling punctures and to establish equality on the frozen charts. We then apply extension of regular functions across subsets of codimension at least two. The complement of the union of cluster tori is studied systematically as the deep locus by Beyer and Muller\cite{beyer25deep}; their framework also applies to intermediate algebras such as skein algebras.

        Therefore, in Section 4, we will establish a lemma for estimating the height of an ideal such as $I_0$. Equivalently, we estimate the (Krull) dimension of the quotient algebra $Sk^{ta}(\Sigma)/I_0$. Due to the equivalence of Krull and Gelfand--Kirillov dimensions, the goal becomes controlling the exponent of growth of some powers of a generating set, which (by applying skein relations) is further reduced to the dimension of curve systems modulo a certain equivalence.\par

		\subsection*{AI disclosure}
No AI tools were used in conducting the research or preparing this manuscript.
        
		\subsection*{Acknowledgments}
		The author would like to extend sincere gratitude to Zhe Sun for introducing him to the subject of skein and cluster algebras and suggesting that he consider the problems in the present paper. The ingenuity of the localized skein and cluster algebras comes from Linhui Shen, Zhe Sun and Daping Weng \cite{shen23SLn} while they were working on the quantization of the $\mathcal{A}$-cluster variety of higher rank. Thang L\^e and Travis Mandel kindly answered questions of the author. The author is indebted to Yunhang Chang, Zhichao Chen and Zhe Sun for enlightening discussions and is grateful to Xiaowu Chen, Hiroaki Karuo and Han-Bom Moon for pertinent advice. Last but not least, the author thanks his beloved parents for consistent support. Any comments and criticisms for the improvement of the present paper are warmly welcomed.
		\section{Preliminaries}

        \subsection{Roger--Yang and tagged skein algebras}

        \begin{defnprop}
        A surface is a connected oriented two-dimensional smooth manifold of finite type, possibly with boundary. Thus $\Sigma$ is diffeomorphic to a compact connected oriented surface $\Sigma_g$ of genus $g$ with finitely many points (called punctures) and finitely many open disks removed; the closures of the disks are pairwise disjoint and contain no punctures. Let a marked surface be a surface with a set $\mathcal{M}$ of finitely many marked points on the boundary $\partial\Sigma$.
        \end{defnprop}

        In our figures, a puncture appears as a hollow circle, while a marked point appears as a solid circle.

        Next, we define one-dimensional objects on a surface following \cite{roger14skein}, such as (simple) arcs, (simple) loops, curve systems, etc. These are special cases of a generalized curve.

        \begin{defn}
        Let $\overline{\Sigma}$ be the compact surface obtained by adding $p$ points at the punctures of $\Sigma$. A generalized curve is an immersion $D\to \overline{\Sigma}$, where $D=\coprod_i [0,1]\amalg\coprod_j S^1$ is a finite collection of closed intervals  and circles, such that: 
        (1) circles and interiors of intervals are mapped to $\Sigma-\partial\Sigma$;\\
        (2) endpoints of intervals are mapped to marked points or punctures.\\
        In a generalized curve, a component corresponding to an interval is called an arc, and a component corresponding to a circle is called a loop. 
        \end{defn}
        
         Generalized curves are unoriented and are considered up to reparametrization of their components.  We denote the number of strands (i.e. number of components when restricted to a small neighborhood) of a generalized curve at a puncture or a marked point by its weight at the puncture or the marked point.
         Next, we introduce an equivalence, regular homotopy, on the set of generalized curves, which is crucial to the definition of skein algebras.

        \begin{defn}
        Two generalized curves are regularly homotopic if there exists a collection $D$ of closed intervals and circles and a map $F:D\times [0,1]\to \overline\Sigma$ such that:\\
        (1) $F(-,0)$ and $F(-,1)$ are the designated generalized curves;\\
        (2) $F$ is a homotopy relative to $\partial D$;\\
        (3) $F(-,t)$ is a generalized curve for every $t\in [0,1]$.
        \label{generalized arc}
        \end{defn}

        The following objects are special cases of a generalized curve. \\
        (1) simple arc: A generalized arc where $D=[0,1]$ is an interval and the map $(0,1)\to \Sigma$ is an embedding;\\
        (2) simple loop: A generalized curve where $D=S^1$ is a circle and the map $S^1\to\Sigma$ is an embedding;\\
        (3) simple diagram: A finite union of simple arcs and simple loops where two components have no common point except at marked points or punctures;\\
        (4) curve system: A simple diagram where no component is contractible and no loop is peripheral (i.e. bounding a once-punctured disk).

        \begin{rem}
        The notion of regular homotopy is not equivalent to homotopy in general (where homotopy is also relative to endpoints, for an arc). Indeed, a simple curve is homotopic to the curve obtained by inserting a small curl, but they are not regularly homotopic. However, if we restrict to the class of simple diagrams, then two simple diagrams are homotopic if and only if they are regularly homotopic. In fact, if we only look at each component (i.e. simple arc or simple loop) that is non-contractible, the notion of homotopy is even equivalent to isotopy, see\cite[Proposition 1.10]{farb12primer}.
        \label{regularhomotopic}
        \end{rem}

        By abuse of notation, we also use a generalized curve to denote its regular homotopy class.

        Following \cite{roger14skein} and \cite{muller16skein}, we recall the classical Roger--Yang skein algebra and its generalization, the classical Muller--Roger--Yang skein algebra $Sk^{RY}_{q\to1}(\Sigma)$. We will call it (Muller--)Roger--Yang skein algebra and denote it by $Sk^{RY}(\Sigma)$ in the following context unless otherwise specified.
        
        \begin{defnprop}[Muller--Roger--Yang skein algebra]
         (\cite{roger14skein}\cite{karuo25cluster}) Let $v_1,...,v_p$ be formal variables, each associated to a puncture. For a surface $\Sigma$, the Muller--Roger--Yang skein algebra $Sk^{RY}(\Sigma)$ is the $\mathbb{C}[v_1^{\pm1},...,v_p^{\pm1}]$-algebra (where the $v_i$'s are formal variables) generated by generalized curves and inverses of boundary arcs (i.e. simple arcs homotopic relative to their endpoints to a boundary segment between consecutive marked points, which may coincide) on $\Sigma$ modulo the skein relations, including the Kauffman bracket skein relation, the puncture skein relations and the framing relations. When $\Sigma$ has no boundary, i.e. when $\Sigma=\Sigma_{g,p}$, the Muller--Roger--Yang skein algebra $Sk^{RY}(\Sigma)$ is also called the Roger--Yang skein algebra. The set of curve systems with weight at most 1 at each puncture forms a basis of $Sk^{RY}(\Sigma)$ over $\mathbb{C}[v_1^{\pm1},...,v_p^{\pm1}]$, with arbitrary integer multiplicities of boundary arcs allowed. Throughout, when curve systems index bases of boundary-localized algebras, negative boundary multiplicities represent inverse powers of boundary arcs.
        \label{RYskein}
        \end{defnprop}

        \begin{overpic}[width=0.5\textwidth]{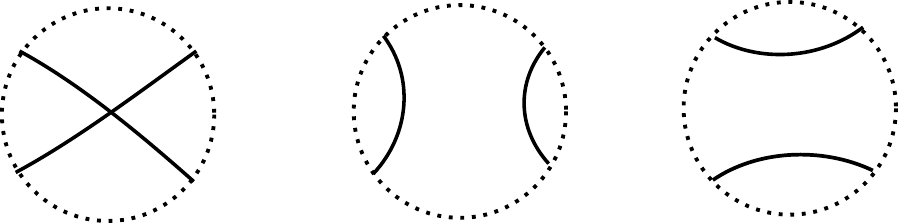}
        \put(110,12){(The Kauffman bracket relation)}
         \put(30, 12){$=$}
         \put(68, 12){$+$}
        \end{overpic}
    
       \begin{overpic}[width=0.5\textwidth]{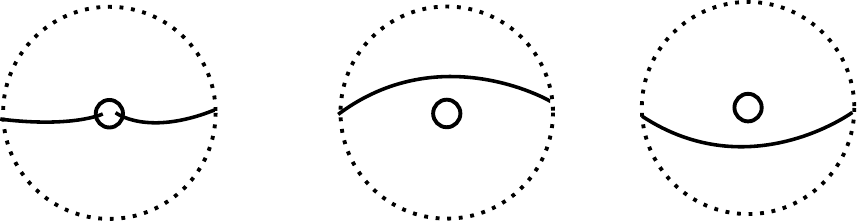}
        \put(110,12){(The puncture skein relation)}
        \put(-7,12){$v_i$}
        \put(30,12){$=$}
        \put(68,12){$+$}
        \end{overpic}

        \begin{overpic}[width=0.35\textwidth]{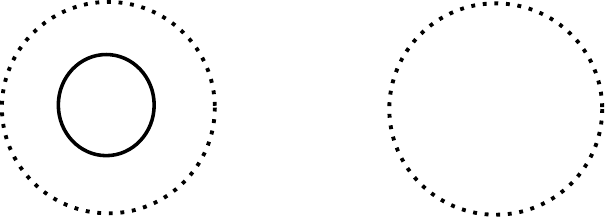}
        \put(157.5,15){(The framing relation)}
        \put(42.5,15){$=$}
        \put(52,15){$-2$}
        \end{overpic}

        \begin{overpic}[width=0.35\textwidth]{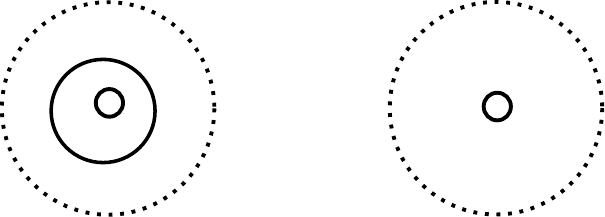}
        \put(157.5,15){(The puncture framing relation)}
        \put(42.5,15){$=$}
        \put(55,15){$2$}
        \end{overpic}

        \begin{overpic}[width=0.12\textwidth]{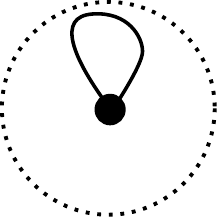}
        \put(460,40){(Vanishing of contractible arc)}
        \put(125,40){$=$}
        \put(160,40){$0$}
        \end{overpic}

        \begin{rem}
        The Muller--Roger--Yang skein algebra in our definition is in fact the boundary-localized Muller--Roger--Yang skein algebra, to be distinguished from the Muller--Roger--Yang skein algebra defined in \cite{bloomquist23degenerations}\cite{karuo25cluster} that does not contain inverses of boundary arcs.
        \end{rem}

        Next, we recall tagged arcs and tagged triangulations introduced by Fomin, Shapiro and Thurston \cite{fomin08cluster1}. These are more natural objects in terms of the compatibility of skein and cluster algebras.

        \begin{defn}[Tagged arcs, compatibility and tagged curve systems]
        A generalized tagged arc is an arc with a decoration, namely plain or notched, at all its ends that point to a puncture. Ends at boundary marked points are always plain.
           A tagged arc is a generalized tagged arc, which is simple, not contractible, not bounding a once-punctured monogon, and has the same decorations at its ends if the two ends are at the same puncture. Say two tagged arcs $\alpha$ and $\beta$ are compatible if their underlying arcs have representatives with disjoint interiors and the following conditions are satisfied:\\
           (1) if the underlying arcs are non-homotopic, the tags agree at every common endpoint;\\
           (2) if the underlying arcs are homotopic, at least one pair of corresponding ends has the same tag.\\
           A tagged curve system is a curve system together with a choice of decoration (plain or notched) at each puncture with positive weight.
        We say two tagged curve systems are (regularly) homotopic if their untagged versions are (regularly) homotopic and the decorations are identical. 
        \label{tagged}
        \end{defn}
      In our figures, a plain arc still appears as the arc itself, while a tagged arc appears as an arc with a marking at each of its ends.

        \begin{rem}
        One may also define a tagged curve system as a collection of tagged arcs and noncontractible nonperipheral simple loops with no interior intersection, such that tagged arcs are pairwise compatible. If one is allowed to apply the digon relation (see Definition-Proposition \ref{taggedskein}), then tagged curve systems in this definition are in bijection with those in the previous one. There is an analogous definition for tagged bracelets in Definition \ref{bracelettaggedbracelet}.
        \end{rem}

        The concept of triangulation is crucial to the study of skein algebras, as well as cluster algebras that arise from moduli of local systems on surfaces. 

        \begin{defn}[Ideal and tagged triangulations]
        We say that a marked surface $\Sigma$ is a triangulable surface if the following conditions hold: (1) $\Sigma$ has at least 1 puncture or boundary marked point, and each boundary component has at least 1 marked point; \\
        (2) $\Sigma$ does not belong to any of the exceptional types of surface: a disk with $\leq 2$ marked points but no puncture, a disk with 1 marked point and 1 puncture, or a sphere with at most 3 punctures.
        An (ideal) triangulation is a maximal collection of pairwise non-homotopic noncontractible simple arcs on a triangulable surface $\Sigma$, including boundary arcs, such that any two arcs can only intersect at boundary marked points or punctures. A tagged triangulation is a maximal collection of distinct pairwise compatible tagged arcs on $\Sigma$. 
        \label{triangulation}
        \end{defn}

        A triangulation must contain all boundary arcs, and has $6g-6+3p+|\mathcal{M}|+3|\partial\Sigma|$ non-boundary arcs, where $|\partial\Sigma|$ is the number of boundary components on $\Sigma$. For a triangulable surface $\Sigma$, a triangulation always exists and divides $\Sigma$ into (possibly self-folded) triangles (i.e. triangles with two sides identified). Similarly, a tagged triangulation contains all boundary arcs, and the number of non-boundary arcs contained in a tagged triangulation is also $6g-6+3p+|\mathcal{M}|+3|\partial\Sigma|$. Tagged triangulations do not use arcs bounding once-punctured monogons, but their associated ideal triangulations may contain self-folded triangles. The ideal triangulation associated to a tagged triangulation is unique up to isotopy and is obtained by forgetting the decorations and possibly changing some arcs subject to the following rule:\\
        (1) If the underlying simple arc of a tagged arc $\alpha$ is homotopic to that of another tagged arc $\beta$, then replace one of the two arcs with the simple arc bounding a once-punctured monogon containing the two arcs, such that $\alpha$ and $\beta$ are tagged differently at the interior puncture of the monogon;\\
        (2) If the underlying simple arc is unique up to homotopy in the triangulation, then it stays unchanged. 

        \begin{overpic}[width=0.3\textwidth]{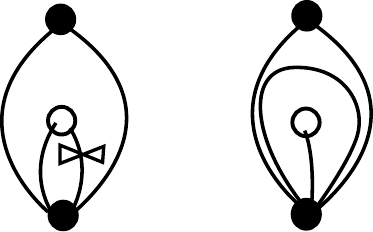}
        \put(120,30){Ideal triangulation and tagged triangulation}
        \put(120,20){of once-punctured disk with 2 marked points}
        \end{overpic}

        To motivate the tagged skein algebra and link it with the Muller--Roger--Yang skein algebra,
        we associate each tagged arc $\tilde\alpha$ with the value $\alpha\prod\limits_{i=1}^p v_i^{k_i}$ if $\alpha$ is the curve system with all decorations forgotten and $k_i$ is the weight of $\alpha$ at the $i$th puncture, with $i$ running over the indices of punctures with notched decorations. Now we make the following definition. Again, we work in the classical setting unless otherwise specified.
       
        \begin{defnprop}
        (\cite{mandel23bracelet}\cite[Theorem 4.1]{musiker13bases}\cite[Theorem 1]{thurston14positive}) For a triangulable marked surface $\Sigma$, the tagged skein algebra $Sk^{ta}(\Sigma)$ is the $\mathbb{C}$-subalgebra of $Sk^{RY}(\Sigma)$ generated by loops, the values of tagged arcs defined above, and inverses of boundary arcs. Equivalently, it is the commutative algebra of generalized tagged diagrams modulo the Kauffman bracket relations, the local digon relations, the (puncture) framing relations and the vanishing of contractible arcs, with boundary arcs inverted; see \cite[Corollary 9.9]{mandel23bracelet}. By \cite[Proposition 3.9]{mandel23bracelet} (or one may simply check the skein relations), all the generalized tagged arcs lie in $Sk^{ta}(\Sigma)$. The collection of tagged curve systems $TCS(\Sigma)$ is a $\mathbb{C}$-linear basis of $Sk^{ta}(\Sigma)$.
        \label{taggedskein}
        \end{defnprop}

        \begin{overpic}[width=0.5\textwidth]{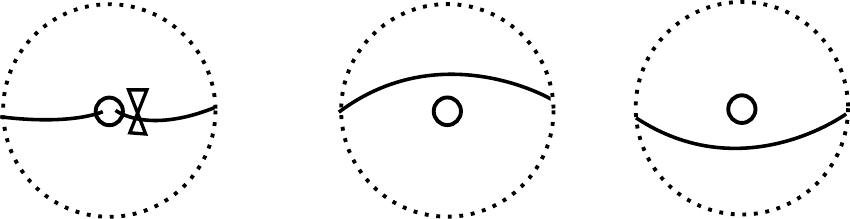}
        \put(110,12){(The digon relation)}
        \put(30,12){$=$}
        \put(68,12){$+$}
        \end{overpic}

        The collection of tagged curve systems forms a basis of the tagged skein algebra, called the bangle basis. In general, this basis need not be positive: the product of two curve systems may have a negative coefficient for some tagged curve system (\cite[Corollary 4.16]{thurston14positive}). To find a positive basis for $Sk^{ta}(\Sigma)$, we invoke the bracelet basis.
        For a noncontractible nonperipheral simple loop $l$ and $k\geq1$, define the bracelet loop $Brac_k(l)$ of weight $k$ to be the standard immersed representative that travels around $l$ $k$ times and has $k-1$ transverse self-intersections, see \cite{musiker13bases}. In the skein algebra, $Brac_k(l)=T_k(l)$, where $T_k$ is the Chebyshev polynomial (which is a unique degree-$k$ monic polynomial with integer coefficients satisfying $T_k(2\cos x)=2\cos{kx}$). 

        \begin{defn}[Bracelet and tagged bracelet]
        A bracelet on $\Sigma$ is a union of simple arcs and bracelet loops whose underlying simple curves have pairwise disjoint interiors, with no parallel bracelet loops, such that no arc or loop is contractible and no loop is peripheral. We allow the weight of a boundary arc in the collection to be negative (so at the level of skein algebra it should be in the denominator if the weight is negative).  A tagged bracelet is a bracelet with all arcs decorated identically at each puncture. 
        \label{bracelettaggedbracelet}
        \end{defn}

        Note that the empty set is also a tagged curve system, and is a tagged bracelet. Both the set of tagged bracelets and the set of tagged curve systems (between which there are simple linear transitions) form a $\mathbb{C}$-linear basis of the tagged skein algebra. The bracelet basis is positive, meaning that its structure constants are nonnegative integers.

        \begin{prop}
         (\cite[Theorem 4.1]{musiker13bases}\cite[Theorem 1]{thurston14positive}\cite{mandel23bracelet}) The collection $TBr(\Sigma)$ of tagged bracelets forms a positive $\mathbb{C}$-linear basis of the tagged skein algebra $Sk^{ta}(\Sigma)$ over $\mathbb{C}$. 
        \label{basis}
        \end{prop}


        To conclude the subsection, we mention some important algebraic properties of the Muller--Roger--Yang skein algebra, which is going to be used in other parts of the paper (e.g. the compatibility of skein and cluster algebras).

        \begin{thm}
        (1) (\cite[Corollary 3]{bloomquist23degenerations}\cite[Theorems 2.8 and 2.11]{moon24consequences}\cite[Theorem A and Proposition 6.4]{karuo25cluster}) For any triangulable surface $\Sigma$, the Muller--Roger--Yang skein algebra $Sk^{RY}(\Sigma)$ and the tagged skein algebra $Sk^{ta}(\Sigma)$ are finitely generated integral domains. \\
        (2) Suppose that $\Sigma$ is not a once-punctured closed surface. Then $Sk^{ta}(\Sigma)$ is a Gorenstein domain and hence satisfies Serre's condition $S_2$: a rational function regular away from a closed subset of codimension at least two is regular globally. This follows from Theorem \ref{ghkkcor}(2)--(4), independently of the equality $Sk^{ta}(\Sigma)=U(\Sigma)$.
        \label{skeinalgebraic}
        \end{thm}

         \begin{rem}
            It is H.-B. Moon who reminded us that the finite generation of $Sk^{ta}(\Sigma)$ can be shown completely skein-theoretically. Indeed, \cite{moon24consequences}\cite{karuo25cluster}\cite{bloomquist23degenerations} show that various kinds of skein algebras are finitely generated, including tagged skein algebras, L\^e--Roger--Yang skein algebras, and Muller--Roger--Yang skein algebras. Later, we shall see that the finite generation of the tagged skein algebra is also implied by that of the middle cluster algebra, introduced by Gross--Hacking--Keel--Kontsevich\cite{gross18canonical}, and the comparison between bracelet and theta bases due to Mandel--Qin\cite{mandel23bracelet}.
        \end{rem}

        \subsection{Cluster algebras and cluster varieties associated to surfaces}
        Let $\Sigma$ be a marked surface, and let $\Sigma'$ be the compact surface with punctures on $\Sigma$ replaced by holes, which is homotopy equivalent to $\Sigma$. Let $T'\Sigma'$ be the punctured tangent bundle of $\Sigma'$, namely the tangent bundle with the zero section removed. A twisted $SL_2$ (viewed as a group scheme over $\mathbb{Q}$)-local system is an $SL_2$-local system on $T'\Sigma'$ such that the monodromy representation of $\mathbb{Z}=\ker (\pi_1(T'\Sigma')\to\pi_1(\Sigma'))$ maps the generator to $-Id\in SL_2$. A spin structure splits the central $\mathbb Z/2\mathbb Z$-extension $\pi_1(T'\Sigma')/\langle z^2\rangle$ of $\pi_1(\Sigma')$, where $z$ is the tangent-fiber generator, and thereby identifies twisted local systems with ordinary local systems. The tangent lift of a contractible embedded circle has monodromy $-I$, whereas the trivial element of the ordinary fundamental group has monodromy $I$. Let $U$ be a maximal unipotent subgroup of $SL_2$, and let $\mathcal F=SL_2/U\simeq\mathbb A^2\setminus\{0\}$ be the principal affine variety. A decoration of a twisted local system $\mathcal L$ is a flat section of the associated principal affine bundle $\mathcal L\times_{SL_2}\mathcal F$ over the boundary collar $\mathcal N=[0,1]\times(\partial\Sigma'-\mathcal M)$, lifted to $T'\Sigma'$ using the oriented boundary tangent direction. Thus a decoration specifies a flat nonzero vector, including its scale, on each boundary interval and a nonzero monodromy-invariant vector around each puncture; see \cite[Definition 2.4]{fock06moduli}.
        \begin{defn}
        Let $\mathcal{A}_{SL_2,\Sigma}$ be the moduli stack of decorated twisted $SL_2$-local systems on $\Sigma$. 
        \end{defn}

        In the absence of boundary marked points, $\mathcal{A}_{SL_2,\Sigma}$ serves as an algebro-geometric avatar of the decorated Teichm\"uller space $\mathcal{T}^d(\Sigma)$ defined by Penner\cite{penner87decorated}, namely the moduli space of marked hyperbolic metrics with finite area together with a choice of length of a horocycle at each puncture.

        \begin{thm}
        (\cite[Theorem 1.7(b)]{fock06moduli}) For a surface $\Sigma=\Sigma_{g,p}$ with $p>0$ and $2g-2+p>0$, fixing a spin structure on $\Sigma$, the decorated Teichm\"uller space is isomorphic to the set of positive real points on the space of decorated twisted $SL_2$-local systems, that is, $$\mathcal{T}^d(\Sigma)=\mathcal{A}_{SL_2,\Sigma}(\mathbb{R}_{>0})$$
        \end{thm}

        Fock and Goncharov showed that the space $\mathcal{A}_{SL_2,\Sigma}$ admits a regular positive atlas. In other words, $\mathcal{A}_{SL_2,\Sigma}$ contains a set of toric charts, called cluster charts, as open subsets that capture the combinatorics of $\mathcal{A}_{SL_2,\Sigma}$. Cluster charts have positive transition functions determined by cluster mutations. Remarkably, the mutation rule of clusters is encoded completely in an initial combinatorial object, which we now formally define.\par 

        \begin{defn}
        Let a quiver $\mathbf{q}$ be the data $(\Lambda,\Lambda_0,\{e_i\},(-,-))$ where $\Lambda$ is a free abelian group of finite rank and $(-,-)$ is a skew-symmetric bilinear form $\Lambda\times\Lambda\to \mathbb{Q}$, $\{e_i\}$ is a basis of $\Lambda$ and a subset of $\{e_i\}$, called the frozen basis vectors (the complement of which consists of the mutable basis vectors), spans $\Lambda_0$. We require that $(e_i,e_j)\in \mathbb{Z}$ unless both $e_i$ and $e_j$ are frozen. The set $I$ parametrizes the  basis vectors $\{e_i\}$ and the set $I_0$ parametrizes the frozen basis vectors.
        \end{defn}

        Let $\epsilon_{ij}=(e_i,e_j)$ and call the matrix $(\epsilon_{ij})_{1\leq i,j\leq n}$ the exchange matrix. If $e_k$ is a nonfrozen basis vector, then the mutation $\mu_k$ in the $k$-direction defines a new quiver 
        $\mathbf{q}'=(\Lambda,\Lambda'_0,\{e'_i\},(-,-))$, where $\Lambda'_0$ is spanned by $\{e'_i:i\in I_0\}$ and $[a]_+=\max(a,0)$, such that \\
        \[
        e_i':=
        \begin{cases} 
            e_i+[\epsilon_{ik}]_{+}e_k & \text{if }  i\neq k \\
            -e_k & \text{if } i=k
        \end{cases}
        \]
        and we calculate the change of the structural constants by \\
        \[
        \epsilon_{ij}':=
        \begin{cases} 
            -\epsilon_{ij} & \text{if }  k\in \{i,j\} \\
            \epsilon_{ij} & \text{if }  k\notin\{i,j\}, \epsilon_{ik}\epsilon_{jk}\geq0 \\
            \epsilon_{ij}-|\epsilon_{ik}|\epsilon_{jk} & \text{if } k\notin\{i,j\}, \epsilon_{ik}\epsilon_{jk}<0
        \end{cases}
        \]

        We write $\mathbf{q}_1\sim\mathbf{q}_2$ if the two quivers can be related by a sequence of mutations.

        Next, we assign each quiver a set of coordinates (namely a collection of $n$ algebraically independent variables over a field $k$ of characteristic 0) say $\{A_i\}_{i\in I}$, called a cluster. Each variable is attached to a distinguished vertex of $\mathbf{q}$. Under the mutation $\mu_k$, the coordinates (also called cluster variables) are transformed according to the following rule: 
        
        \[
        A_i'=
        \begin{cases}
        A_k^{-1}(\prod\limits_{j|\epsilon_{jk}<0} A_j^{-\epsilon_{jk}}+\prod\limits_{j|\epsilon_{jk}>0} A_j^{\epsilon_{jk}}) &\text{ if } i=k\\
        A_i &\text{ if } i\neq k   
        \end{cases}
        \]

        It is important to notice that mutating the same vertex twice is an isomorphism of seeds. Such an operation induces the identity map on cluster variables and the corresponding exchange matrix, although not the identity on the vector space spanned by the basis vectors.

        \begin{defn}
        A seed $\mathbf{s}$ contains the data of a quiver and its cluster. A mutation of a seed is a mutation of its quiver and a replacement of its cluster according to the rule above. 
        \end{defn}

        To encode the combinatorics of seeds under mutations, we recall the exchange graph of a quiver (or its equivalence class). The symmetry of the exchange graph is characterized by the cluster modular group.

        \begin{defn}[Exchange graph and cluster modular group]
        The exchange graph of a quiver $\mathbf{q}$ is the $k$-regular graph (where $k$ is the number of mutable vertices, called the rank of $\mathbf{q}$) whose vertices are the seeds mutation equivalent to an initial seed, considered up to simultaneous relabeling of their mutable variables and exchange matrices, and whose edges correspond to mutations; see \cite[Definition 5.4]{fomin08cluster1}. The cluster modular group consists of compositions of mutations and seed isomorphisms returning the initial quiver to itself, modulo those inducing the identity on both the $\mathcal A$- and $\mathcal X$-tori; see \cite[Section 1.2]{fock09cluster}.
        \end{defn}

        Note that the cluster modular group acts as a group of automorphisms of the function field $\mathbb{C}(A_i)(i\in I)$. \par
        Now we are ready to define the main objects studied in cluster theory, namely the cluster algebra, upper cluster algebra, and a geometric version, called cluster variety.
      
        The transformation rule between clusters of two adjacent seeds  defines a birational map $\mathrm{Spec}(\mathbb{C}[A_i^{\pm1}])\dashrightarrow\mathrm{Spec}(\mathbb{C}[A_i'^{\pm1}])$ between two tori; such a torus is called a cluster chart and is denoted by $T_\mathbf{s}$. (Despite its notation, a cluster chart depends only on the cluster rather than on the whole seed.) For any two (not necessarily adjacent) cluster charts we have such a birational map, and we can take the birational map that is an isomorphism on the largest possible open sets. If we glue all the tori via these birational maps, according to \cite[Proposition 2.4]{gross15birational}, we would obtain a scheme $\mathcal{A}$, called the $\mathcal{A}$-cluster variety (or simply cluster variety unless otherwise specified) with respect to the initial quiver. 
        Let the (ordinary) cluster algebra $A$ be the subalgebra of $\mathbb{C}(A_i)$ ($i\in I$) generated over $\mathbb{C}$ by all the cluster variables and the inverses of frozen variables, and the upper cluster algebra $U=\Gamma(\mathcal{A},\mathcal{O}_{\mathcal{A}})$ be the ring of global sections of the scheme $\mathcal{A}$. By definition, we also have $$U=\bigcap\limits_{\mathbf{s}'\sim\mathbf{s}}\mathbb{C}[A_i'^{\pm1}](i\in I),$$ that is, the rational functions in $\mathbb{C}(A_1,...,A_n)$ which are Laurent polynomials when expressed in any set of coordinate functions. Furthermore, by the Laurent phenomenon\cite{fomin01cluster}\cite[Corollary 3.11]{gross15birational} we have $A\subseteq U$. 
        It is clear that mutation-equivalent quivers give rise to the same $\mathcal{A}$-cluster variety, cluster algebra and upper cluster algebra. 

        Now it comes to relating a triangulation of a triangulable surface $\Sigma$ to a quiver. The edges in a triangulation correspond to the vertices of the quiver $\mathbf{q}$, and the exchange matrix of the quiver should be exactly the adjacency matrix of the triangulation (to be defined below). The operation of mutating a vertex should correspond to flipping an edge (i.e. the operation of replacing a diagonal of a quadrilateral by the other diagonal). However, when $\Sigma$ has punctures, there may appear self-folded triangles in an ideal triangulation, and the internal edge cannot be flipped. This is why we work with tagged triangulations instead of ideal triangulations. \par
        A flip at a non-boundary tagged arc $\alpha$ of a tagged triangulation $T$ replaces it by the unique different tagged arc $\alpha^*$ for which $(T\setminus\{\alpha\})\cup\{\alpha^*\}$ is a tagged triangulation; see \cite[Theorem 7.9]{fomin08cluster1}. For an ordinary flippable edge, this is the usual exchange of diagonals in a quadrilateral.\par
         
        To each tagged triangulation, we attach a quiver $\mathbf{q}$ whose mutable vertices correspond to tagged arcs and whose frozen vertices correspond to boundary arcs, and the exchange matrix $(\epsilon_{ij})_{1\leq i,j\leq n}$ equals the adjacency $B$ matrix of $\Sigma$, where the adjacency matrix of $\Sigma$ with respect to a tagged triangulation is defined to be that with respect to the associated ideal triangulation. We recover an ideal triangulation (see Definition \ref{triangulation} and below) from the tagged triangulation and number the arcs in the tagged triangulation by $\alpha_1,...,\alpha_n$, the corresponding arcs in the ideal triangulation by $\beta_1,...,\beta_n$. Then, if $\beta_i$ is a self-folded edge then let $\pi_T(\beta_i)$ be the arc bounding the monogon, otherwise set $\pi_T(\beta_i)=\beta_i$, and denote the collection of unfolded triangles in the associated ideal triangulation by $T$, then $b_{ij}=\sum\delta_{ij,\Delta}$ where the summation runs over all unfolded triangles of $T$, and 
        \[       
         \delta_{ij,\Delta}=
        \begin{cases}
        0 & \text{if } \pi_T(\beta_i)\notin\Delta \text{ or } \pi_T(\beta_j)\notin\Delta \text{ or } \pi_T(\beta_i)=\pi_T(\beta_j)\\
        -1 & \text{if } \pi_T(\beta_i),\pi_T(\beta_j)\in\Delta \text{ and } \pi_T(\beta_i) \text{ is clockwise to } \pi_T(\beta_j) \text{ in } \Delta\\
        1 & \text{if } \pi_T(\beta_i),\pi_T(\beta_j)\in\Delta \text{, }\pi_T(\beta_i)\text{ counterclockwise to }\pi_T(\beta_j)
        \end{cases}
        \]

        We define a quiver associated to a triangulation as that with exchange matrix $\epsilon_{ij}=b_{ij},1\leq i,j\leq n$. It turns out that the flip of a tagged triangulation corresponds to the mutation of a quiver. Formally, we have

        \begin{prop}
        (\cite[Proposition 4.8, Lemma 9.7]{fomin08cluster1})
        Let $\mathbf{q}$ be the quiver associated to a tagged triangulation of $\Sigma$, then the mutation $\mu_k(\mathbf{q})$ at a mutable vertex $k$ has the same exchange matrix as the quiver associated to the tagged triangulation with its $k$th arc flipped.
        \label{mutationflip}
        \end{prop}

        \begin{rem}
        Tagged arcs ensure that every non-boundary arc can be flipped. The entries of the adjacency matrix are computed from the associated ideal triangulation by the signed-angle rule above, including the map $\pi_T$ for self-folded triangles.
        \end{rem}

        Moreover, for most triangulable surfaces, the complex of tagged arcs (with two tagged arcs adjacent if they are compatible) turns out to be isomorphic to the complex of cluster variables (with two cluster variables adjacent if they are contained in the same cluster) under the identification of cluster variables as tagged arcs. These complexes control the combinatorics of triangulations (resp. clusters) and flips (resp. mutations).

        \begin{thm}
        (\cite[Theorem 5.6]{fomin08cluster1}\cite[Theorem 7.11]{fomin08cluster1}\cite[Corollary 6.2]{fomin08cluster2}) When $\Sigma\neq\Sigma_{g,1}$, the cluster complex is isomorphic to the tagged arc complex, which is connected; when $\Sigma=\Sigma_{g,1}$, the cluster complex is isomorphic to one of the two isomorphic connected components of the tagged arc complex, namely the complex of plain arcs and the complex of notched arcs.
        \label{arcvariable}
        \end{thm}

        We denote the $\mathcal{A}$-cluster variety associated to a quiver of a tagged triangulation of $\Sigma$ by $\mathcal{A}$ unless otherwise specified. The corresponding cluster algebra and the upper cluster algebra are denoted by $A(\Sigma)$ and $U(\Sigma)$. In the following context, we shall take the ground field $k$ to be $\mathbb{C}$ unless otherwise specified and base change $\mathcal{A}_{G,\Sigma}$ from $\mathbb{Q}$ to $\mathbb{C}$. The relation with the moduli space $\mathcal{A}_{SL_2,\Sigma}$ is provided by its positive atlas.

        \begin{thm}
        (\cite[Theorem 1.5]{fock06moduli}) The moduli space $\mathcal{A}_{SL_2,\Sigma}$ has a regular positive atlas: its ideal-triangulation coordinate charts are open torus embeddings, with positive birational transition maps.
        \end{thm}

        In what follows, $\mathcal A$ denotes the cluster variety defined by gluing seed tori. The local-system charts provide its geometric interpretation; we do not identify this cluster variety with the entire moduli stack $\mathcal{A}_{SL_2,\Sigma}$.\par

        In \cite{muller13acyclic}, Muller studied a special kind of cluster algebra, called a locally acyclic cluster algebra, which incorporates the class of acyclic cluster algebras (cluster algebras that come from a quiver with no oriented cycle in the mutable part), and the class of isolated cluster algebras (ones that come from a quiver with no arrows between mutable vertices). A locally acyclic cluster algebra is a cluster algebra whose spectrum admits a finite open cover by spectra of acyclic (equivalently, isolated \cite[Theorem 6.5]{muller13acyclic}) cluster localizations. Here a cluster localization is obtained by freezing a subset of the variables of a seed and agrees with the localization inverting those variables. All frozen variables remain inverted, as in our convention above. As an example, we have:

        \begin{prop}
(1) (\cite[Theorem 10.6]{muller13acyclic}) If $\Sigma$ has at least two boundary marked points, its cluster algebra is locally acyclic. The same conclusion holds for an ice quiver with the same mutable part, with frozen variables inverted.

(2) (\cite[Proposition 12 and its proof]{canakci15cluster}) If $\Sigma$ has exactly one boundary marked point and at least one puncture, its mutable quiver is Louise. Thus any ice quiver with this mutable part, with frozen variables inverted, defines a locally acyclic cluster algebra.
\label{surfacelocallyacyclic}
\end{prop}

        Locally acyclic cluster algebras enjoy many nice algebraic properties. For instance,

        \begin{thm}
        (1)(\cite[Theorem 2]{muller14au}) For a locally acyclic cluster algebra, $A=U=\Gamma(\mathcal{A},\mathcal{O}_{\mathcal{A}})$. That is, the cluster algebra coincides with the upper cluster algebra.

        (2)(\cite[Theorem 4.2]{muller13acyclic}) If $A$ is a locally acyclic cluster algebra, then it is finitely generated, has only complete intersection singularities and is normal (hence is $S_2$).
        \label{locallyacyclicau}
        \end{thm}

        Besides the case of locally acyclic cluster algebras, there is another class of cluster algebras of surface type where the cluster and upper cluster algebras coincide.
        \begin{prop}
(\cite[Theorem 5]{canakci15cluster}) Let $\Sigma$ be an unpunctured triangulable surface with exactly one marked point. For the cluster algebra with trivial coefficients, $A(\Sigma)=U(\Sigma)$.
\label{1markedpointau}
\end{prop}

        \subsection{Freezing, cutting and localization}
        For a quiver $\mathbf{q}$, let $I$ be the set of vertices and $I_0$ be the subset of frozen vertices. Given a subset of mutable vertices $I_1\subseteq I-I_0$, we construct a quiver $\mathbf{q}^{I_1}$ with the same basis and skew-symmetric bilinear form but with frozen vertices $I_1\cup I_0$. This is called freezing $\mathbf{q}$ at $I_1$. Starting from the same cluster, the allowed mutation sequences avoid $I_1$, and the newly frozen variables are inverted. For the surface seeds considered here, the resulting exchange graph is the subgraph obtained by keeping these variables fixed.

        When the cluster algebra is associated to a triangulable surface $\Sigma$, the notion of freezing turns out to be closely related to the cutting of $\Sigma$. Given a non-boundary simple arc $\alpha$ which is neither contractible nor the boundary of a once-punctured disk, cutting $\Sigma$ along $\alpha$ and adding marked points to the preimage of punctures at which $\alpha$ ends, we get a new marked surface (or a disjoint union of two marked surfaces) $Cut_{\alpha}(\Sigma)$, together with a gluing map $f:Cut_{\alpha}(\Sigma)\to \overline{\Sigma}$, where punctures are restored in the target. By assumption, each connected component of $Cut_{\alpha}(\Sigma)$ is triangulable (which one may check by hand), and $f$ is a diffeomorphism away from $f^{-1}(\alpha)$, while it is a double cover along the interior of $f^{-1}(\alpha)$ (which becomes a part of $\partial Cut_{\alpha}(\Sigma)$). It is important to note that, under our definition of triangulable surface, each component of $Cut_{\alpha}(\Sigma)$ is a triangulable surface whenever $\Sigma$ is a triangulable surface.

        If $\alpha$ is the ordinary version of a tagged arc $\tilde{\alpha}$, let $T$ be a tagged triangulation containing $\tilde{\alpha}$. First change all tags at an endpoint puncture whenever necessary to make $\tilde{\alpha}$ plain at both ends. If $T$ contains its conjugate tagged arc, replace that conjugate by the loop enclosing $\alpha$ and the puncture where their tags differ. Taking the preimage under $f$ then gives a tagged triangulation $T'$ of $Cut_{\alpha}(\Sigma)$; endpoints on the newly created boundary are plain. The two copies of $\alpha$ are boundary arcs of the cut surface.

        We may consider the quiver $\mathbf{q}$ associated to a tagged triangulation $T$ (satisfying the above assumptions) of $\Sigma$ and the corresponding quiver $\mathbf{q}'$ of $T'$ on $Cut_{\alpha}(\Sigma)$ (where we may trivially extend our construction of a quiver to multiple components). The quiver $\mathbf{q}'$ is closely related to the freezing $\mathbf{q}^i$ of $\mathbf{q}$ at the vertex $i$ of $\tilde{\alpha}$.

        \begin{prop}
        For a quiver $\mathbf{q}$ that comes from a triangulated surface $\Sigma$, suppose there is a tagged arc $\alpha$ (whose corresponding vertex in $\mathbf{q}$ is $i$) satisfying the assumptions above:
        (1) For $\mathbf{q}_1$ and $\mathbf{q}_1'$ which are mutation-equivalent to the frozen quiver $\mathbf{q}^{i}$ and $\mathbf{q}'$ respectively by the same sequence of mutations avoiding $i$, the mutable parts $\mathbf{q}_{1mut}$ and $\mathbf{q}_{1mut}'$ are isomorphic. Moreover, $\mathbf{q}_1$ and $\mathbf{q}_1'$ have the same exchange graph.\\
        (2) There are vertices $i_1$ and $i_2$ of $\mathbf{q}'$ corresponding to $f^{-1}(\alpha)$, such that the frozen vertices of $\mathbf{q}'$ are exactly the frozen vertices of $\mathbf{q}^i$ with $i$ replaced by $i_1$ and $i_2$. 
        \label{freezingcuttingequality}
        \end{prop}

        \begin{proof}
        The second statement of (1) follows from the identification of a cluster with its triangulation, and the fact that tagged triangulations of $Cut_{\alpha}(\Sigma)$ and tagged triangulations of $\Sigma$ with $\alpha$ frozen are in canonical correspondence. Together with Proposition \ref{mutationflip}, the first statement of (1) follows. (2) follows because the two copies of $\alpha$ become boundary arcs, while the original boundary arcs remain unchanged.
        \end{proof}

        The operation of freezing is also related to localization of cluster and upper cluster algebras.

        \begin{prop}
        Let $\mathbf{s}$ be a seed with quiver $\mathbf{q}$, and let $I$ be a subset of its unfrozen vertices. Let $A^{I}$ and $U^{I}$ be the cluster and upper cluster algebras of the seed obtained by freezing $I$. Let $S=\{A_i:i\in I\}$ be the corresponding variables in $\mathbf{s}$. Then $A[x^{-1}]_{x\in S}$ and $U[x^{-1}]_{x\in S}$ are the localizations of $A$ and $U$ at the multiplicative subset generated by $S$, respectively. There are canonical inclusions $$A^I\subseteq A[x^{-1}]_{x\in S}\subseteq U[x^{-1}]_{x\in S}\subseteq U^I.$$
        \end{prop}

        \begin{proof}
        This is by the definition of cluster and upper cluster algebras, or one may consult \cite[Proposition 3.1]{muller13acyclic} for a proof.
        \end{proof}

         In particular, we have the following corollary.
         \begin{cor}
         In the notation of the last proposition, if $A^I=U^I$, then the four algebras, $A^I$, $A[x^{-1}]_{x\in S}$, $U[x^{-1}]_{x\in S}$, $U^I$ are equal.
         \label{freezinglocalizationcoincide}
         \end{cor}

         In fact, at least in some special cases, the freezing operation is not only defined for a vertex of a specific quiver. Rather, it is defined for a cluster variable regardless of the choice of quiver in a mutation equivalence class.

         \begin{prop}
         Let a quiver $\mathbf{q}$ be associated to a marked surface, and let $\mathbf{q}'$ be a mutation-equivalent quiver. If the seeds of $\mathbf{q}$ and $\mathbf{q}'$ both have a cluster variable $x$, then the frozen quivers at the corresponding vertex are mutation equivalent (up to a change of indices). In particular, the frozen quivers have the same $\mathcal{A}$-cluster variety, cluster algebra and upper cluster algebra. 
         \end{prop}

         \begin{proof}
         Since tagged arcs (when $\Sigma=\Sigma_{g,1}$, consider plain arcs) are in correspondence with cluster variables, we have two tagged triangulations of $\Sigma$ which correspond to $\mathbf{q}$ and $\mathbf{q}'$, respectively. The two triangulations share the tagged arc corresponding to $x$; let $\alpha$ be its underlying ordinary arc. Since each component of $Cut_{\alpha}(\Sigma)$ is a triangulable surface which is not a $\Sigma_{g,1}$, the proposition follows from the fact that any two tagged triangulations of $Cut_{\alpha}(\Sigma)$ (which correspond to tagged triangulations of $\Sigma$ containing $\alpha$) can be related by a sequence of flips.
         \end{proof}

        \subsection{Fock--Goncharov duality conjecture and compatibility of skein and cluster algebras}

        An $\mathcal{A}$-cluster variety is associated with a dual, called the $\mathcal{X}$-cluster variety, defined as a union of toric charts subject to a dual transformation rule. Given a seed, we assign a set of algebraically independent variables $\{X_i\}_{1\leq i\leq n}$ called the $\mathcal{X}$-coordinates, to be distinguished from $\{A_i\}_{1\leq i\leq n}$. The transition rule of $\mathcal{X}$-coordinates under cluster mutation $\mu_k$ is given by
        \[
            X_i'=
            \begin{cases}
            X_k^{-1} &\text{ if } i=k \\
            X_i(1+X_k^{-sgn(\epsilon_{ik})})^{-\epsilon_{ik}} &\text{ if } i\neq k
            \end{cases}
        \]
        We denote the cluster variety obtained by gluing toric charts \[\mathrm{Spec}(\mathbb{C}[(X_1')^{\pm1},\ldots,(X_n')^{\pm1}])\] associated to all seeds by $\mathcal{X}$.\par 
        Following\cite[Section 4]{fock06moduli}\cite[Section 2]{berenstein96parametrizations} a semifield is a set $P$ endowed with two operations, called addition and multiplication. Addition is commutative and associative, multiplication makes $P$ an abelian group, and multiplication distributes over addition. For instance, the sets $\mathbb{Z}^t,\mathbb{Q}^t,\mathbb{R}^t$ of integers, rational numbers and real numbers are semifields under the operations $(\max,+)$; also, $\mathbb{Q}_{>0}$ and $\mathbb{R}_{>0}$ are semifields under the usual addition and multiplication $(+,\cdot)$.\par

        Since the transition functions for $\{A_i\}_{1\leq i\leq n}$ and $\{X_i\}_{1\leq i\leq n}$ are subtraction-free (or positive), the transition rule remains valid for any tropical semifield if we replace $(+,\cdot)$ by the addition and multiplication of the semifield, see\cite[Section 4]{fock06moduli}. Analogously, the spaces $\mathcal{A}$ and $\mathcal{X}$ taking values in $\mathbb{Z}^t$ are denoted by $\mathcal{A}(\mathbb{Z}^t)$ and $\mathcal{X}(\mathbb{Z}^t)$. There are similar notations for $\mathbb{Q}^t$ and $\mathbb{R}^t$.

        In \cite[Section 12]{fock06moduli}, Fock and Goncharov construct a map $\mathbb{I}_{\mathcal{X}}:\mathcal{X}(\mathbb{Z}^t)\to \Gamma(\mathcal{A},\mathcal{O}_{\mathcal{A}})$ where $\Sigma=\Sigma_{g,p}$ and $\mathcal X,\mathcal A$ are the cluster varieties associated to $\mathcal{X}_{PGL_2,\Sigma}$ and $\mathcal{A}_{SL_2,\Sigma}$, respectively. $\mathcal{X}(\mathbb{Z}^t)$ parametrizes $\mathcal{X}$-laminations on $\Sigma$, namely a curve system with positive integral component weights and an orientation at each puncture incident to an arc. Denote this set by $XL(\Sigma,\mathbb Z)$. Combined with the definition of a tagged bracelet, there is a bijection $$\mathcal{X}_{PGL_2,\Sigma}(\mathbb{Z}^t)\cong XL(\Sigma,\mathbb Z)\cong TBr(\Sigma)$$ between the set of integral $\mathcal{X}$-laminations and the set of tagged bracelets, giving a bracelet a plain tag at a puncture if the orientation at the puncture coincides with that of $\Sigma$, and a notched tag otherwise; moreover, a weight-$k$ simple loop (denote the simple loop by $l$) is identified with the bracelet loop $Brac_k(l)$.
        
        We use tangent lifts for loops and good lifts for arcs, bending endpoint tangents to the positive boundary tangent as in \cite[Section 8.1]{costantino22stated}. With this sign convention, the Fock--Goncharov map is characterized by the following axioms. \\
        (1) For a bracelet loop $l$ which is $k$-fold a simple loop $l_0$, $\mathbb{I}_{\mathcal{X}}(l)$ is defined to be the trace of the $k$th power of the monodromy of the tangent lift of $l_0$, namely $\mathbb{I}_{\mathcal{X}}(l)=\lambda^k+\lambda^{-k}$ where $\mathbb{I}_{\mathcal{X}}(l_0)=\lambda+\lambda^{-1}$.\\
        (2) For a simple arc $l$ from puncture $1$ to puncture $2$ (which may coincide), plain at both ends, let $\mathbf w_1$ and $\mathbf w_2$ be the nonzero vectors specified by the decorations. Then $\mathbb I_{\mathcal X}(l)=\omega(\mathbf w_2,\rho(\widehat l)\mathbf w_1)$, where $\widehat l$ is the good lift and $\omega$ is the invariant volume form. The identity $\rho(\widehat{l^{-1}})=-\rho(\widehat l)^{-1}$ makes this evaluation independent of the chosen orientation of $l$. \\
        (3) For two disjoint tagged bracelets $\alpha_1$ and $\alpha_2$ such that $\alpha_1\cup\alpha_2$ is still a tagged bracelet, $\mathbb{I}_{\mathcal{X}}(\alpha_1\cup\alpha_2)=\mathbb{I}_{\mathcal{X}}(\alpha_1)\mathbb{I}_{\mathcal{X}}(\alpha_2)$\\
        (4) The group $(\mathbb{Z}/2\mathbb{Z})^p$ acts on $\mathcal{A}$ birationally by changing the decoration at the $i$th ($1\leq i\leq p$) puncture. More precisely, in a volume-one basis $(\mathbf u_i,\mathbf w_i)$ whose second vector is the decoration, the monodromy of the tangent lift of the $i$th oriented boundary circle is 
        \[
        \begin{bmatrix}
             1 & 0\\
             v_i & 1\\
        \end{bmatrix}
        \]
        In an ideal triangulation, one computes $v_i=\sum_c A_{a_c}/(A_{b_c}A_{d_c})$, where $c$ runs over the triangle corners at the $i$th puncture, $a_c$ is the opposite side, and $b_c,d_c$ are the two incident sides. Corners and sides are counted with multiplicity, and $A_e$ denotes the evaluation in (2) of the side $e$. The action of changing the decoration at the $i$th puncture is defined by $\mathbf{w}_i\mapsto v_i\mathbf{w}_i$, where $\mathbf w_i$ is the nonzero monodromy-invariant vector specified by the decoration at the $i$th puncture. The action multiplies an arc which is plain at the puncture of weight $k$ by $v_i^k$, and takes $v_i$ to $v_i^{-1}$. It is not difficult to check that this birational action is well-defined and is an involution, and such $\mathbb{Z}/2\mathbb{Z}$-actions at different punctures commute, hence a $(\mathbb{Z}/2\mathbb{Z})^p$ action is defined. If a tagged bracelet $l$ is notched at some puncture, one defines $\mathbb{I}(l)$ as the function obtained by applying the action of $(\mathbb{Z}/2\mathbb{Z})^p$ to $\mathbb{I}(l')$, where $l'$ is the plain version of $l$. 

        \begin{rem}
        The birational $(\mathbb{Z}/2\mathbb{Z})^p$ action mentioned above is the Weyl group action on the moduli of local systems. When $\Sigma\neq\Sigma_{g,1}$, tag changes preserve the cluster atlas and act regularly on $\mathcal A$. For $\Sigma=\Sigma_{g,1}$, tag change exchanges the plain and notched triangulation components and is not a cluster transformation of the ordinary atlas. Nevertheless, the evaluations of tagged arcs and loops are Laurent polynomials in every ordinary cluster chart, so every tagged bracelet defines a regular function on $\mathcal A$ \cite[Proposition 3.18]{mandel23bracelet}.

          The construction is natural under orientation-preserving mapping classes and equivariant under the birational tag changes by axiom (4). Thus it is equivariant under the actions generated by these transformations.
        \end{rem}
        
        Since the map $\mathbb{I}_{\mathcal{X}}$ preserves the skein relations (see \cite{fock06moduli}\cite{mandel23bracelet}), it induces a map $\nu:Sk^{ta}(\Sigma)\to \Gamma(\mathcal{A},\mathcal{O}_{\mathcal{A}})=U(\Sigma)$. The skein--cluster comparison shows that $\nu$ is injective \cite[Proposition 3.18]{mandel23bracelet}\cite{moon24consequences}. Since all cluster variables in $U(\Sigma)$ are represented by tagged arcs, the ordinary cluster algebra $A(\Sigma)$ is identified as a subalgebra of the tagged skein algebra.

        \begin{prop}[Compatibility of skein and cluster algebras for $\Sigma_{g,p}$]
        (\cite[Proposition 3.18]{mandel23bracelet}\cite{moon24consequences}) For $\Sigma=\Sigma_{g,p}$, there are inclusions $A(\Sigma)\subseteq Sk^{ta}(\Sigma)\subseteq U(\Sigma)$, identifying cluster variables with tagged arcs if $\Sigma\neq\Sigma_{g,1}$. When $\Sigma=\Sigma_{g,1}$, the cluster variables are identified as the set of plain (or notched) arcs.
        \label{inclusions}
        \end{prop}

       Let $U(\Sigma)^{+}$ be the subset of $U(\Sigma)$ consisting of Laurent polynomials with nonnegative integral coefficients in every cluster chart. Fock and Goncharov conjectured that the set $\mathcal{X}(\mathbb{Z}^t)$ parametrizes the set $\mathbf{L}(\mathcal{A})$ of extremal elements in $U(\Sigma)^{+}$, namely the nonzero elements that cannot be written as a sum of two nonzero elements in $U(\Sigma)^{+}$.  

       \begin{conj}[Fock--Goncharov duality conjecture for $\mathcal{A}_{SL_2,\Sigma_{g,p}}$]
        (1) (\cite[Conjecture 12.2]{fock06moduli}) The image of tagged bracelets is precisely the set of extremal elements in $U(\Sigma)^{+}$. In particular, every element of $U(\Sigma)^{+}$ is a nonnegative integral linear combination of these images.\\
        (2) (\cite[Section 12, Part 4]{fock06moduli}) $\nu:Sk^{ta}(\Sigma)\to U(\Sigma)$ is surjective.
       \label{fgdualityconj}
       \end{conj}


       \begin{rem}
        The conjecture fails for the once-punctured torus $\Sigma_{1,1}$ \cite{zhou20cluster}\cite[Remark 3.19]{mandel23bracelet}. An enlarged atlas is considered in \cite{zhou20cluster}; the two tagged components for general $\Sigma_{g,1}$ are discussed in Section 7 of the present paper.
       \end{rem}

       More generally, for a connected, simply connected, split semisimple algebraic group $G$ over $\mathbb{Q}$, let $^{L}G$ be its Langlands dual group. Let $\mathcal{A}_{G,\Sigma}$ be the moduli space of decorated twisted $G$-local systems on $\Sigma$ and $\mathcal{X}_{^{L}G,\Sigma}$ be the moduli space of framed $^{L}G$-local systems on $\Sigma$. Fock and Goncharov\cite[Section 12]{fock06moduli} conjecture a map $\mathbb{I}_{\mathcal{A}}:\mathcal{A}_{G,\Sigma}(\mathbb{Z}^t)\to \Gamma(\mathcal{X}_{^{L}G,\Sigma},\mathcal{O}_{\mathcal{X}_{^{L}G,\Sigma}})$ and $\mathbb{I}_{\mathcal{X}}:\mathcal{X}_{^{L}G,\Sigma}(\mathbb{Z}^t)\to \Gamma(\mathcal{A}_{G,\Sigma},\mathcal{O}_{\mathcal{A}_{G,\Sigma}})$ satisfying properties similar to those mentioned in Conjecture \ref{fgdualityconj}. The assertion for $\mathbb I_{\mathcal X}$ in this Fock--Goncharov duality conjecture reduces to Conjecture \ref{fgdualityconj} in the case of $G=SL_2$ (i.e. $^{L}G=PGL_2$) and $\Sigma=\Sigma_{g,p}$.

       When $\Sigma$ is replaced by an arbitrary triangulated surface, there is an analogous map $Sk^{ta}(\Sigma)\to U(\Sigma)=\Gamma(\mathcal{A},\mathcal{O}_{\mathcal{A}})$ that can be constructed in the same spirit, whose image contains the ordinary cluster algebra. The cluster algebra $A(\Sigma)$ admits a homomorphism to the tagged skein algebra $Sk^{ta}(\Sigma)$ identifying a cluster variable with the corresponding tagged arc, since the cluster exchange relation is encoded in the Kauffman bracket relation and the digon relation (see \cite{karuo25cluster} for details). The tagged skein algebra $Sk^{ta}(\Sigma)$ admits a homomorphism $\nu:Sk^{ta}(\Sigma)\to U(\Sigma)$, which is related to the fact that the lambda function of a tagged arc or loop is a Laurent polynomial with respect to any triangulation (see \cite{roger14skein}). Intuitively, given a tagged triangulation $T$ and a tagged arc $\alpha$, one reduces the number of interior intersections (between $\alpha$ and $T$) and incompatible decorations by multiplying $\alpha$ by arcs in $T$ sufficiently many times until all skein resolutions are monomials in $T$. The comparison theorem of \cite{karuo25cluster} shows that the maps $A(\Sigma)\to Sk^{ta}(\Sigma)\to U(\Sigma)$ are inclusions. In summary, we have the following proposition.

       \begin{prop}[Compatibility of skein and cluster algebras for triangulable $\Sigma$]
       (\cite[Theorem B (2)]{karuo25cluster}) For any triangulable surface $\Sigma$ with boundary, there is a canonical homomorphism $\nu:Sk^{ta}(\Sigma)\to U(\Sigma)$ identifying tagged arcs with cluster variables and boundary arcs with frozen variables. The map $\nu$ is injective and there are inclusions $A(\Sigma)\subseteq Sk^{ta}(\Sigma)\subseteq U(\Sigma)$, where $Sk^{ta}(\Sigma)$ is identified with its image under $\nu$.
       \label{borderedinclusions}
       \end{prop}

       In the present paper, we give a geometric version of the inclusion of algebras $Sk^{ta}(\Sigma)\to U(\Sigma)$. We recall the construction of $\mathcal{A}=\bigcup_{\mathbf{s}} T_{\Sigma,\mathbf{s}}$ as a union of cluster charts. For each cluster chart $\{A_1,...,A_n\}$, since $\mathbb{C}[A_1,...,A_n]\subseteq Sk^{ta}(\Sigma)\subseteq \mathbb{C}[A_1^{\pm1},...,A_n^{\pm1}]$, it turns out that \[\mathbb{C}[A_1^{\pm1},\ldots,A_n^{\pm1}]=Sk^{ta}(\Sigma)[A_1^{-1},\ldots,A_n^{-1}]\] is a localization of $Sk^{ta}(\Sigma)$. In other words, there is an open immersion $T_{\Sigma,\mathbf{s}}\to \mathrm{Spec}(Sk^{ta}(\Sigma))$ from the torus to the spectrum of the skein algebra. For any two cluster charts $T_{\Sigma,\mathbf{s}_1}$, $T_{\Sigma,\mathbf{s}_2}$, the birational map determined by cluster transformation (which is defined on the largest open subsets as a regular map) is the same as the birational map induced by the open immersions $T_{\Sigma,\mathbf{s}_1}\to \mathrm{Spec}(Sk^{ta}(\Sigma))$ and $T_{\Sigma,\mathbf{s}_2}\to \mathrm{Spec}(Sk^{ta}(\Sigma))$. This says that the open immersions glue into an open immersion $\mathcal{A}\to \mathrm{Spec}(Sk^{ta}(\Sigma))$. One checks easily that the induced map $Sk^{ta}(\Sigma)\to U(\Sigma)=\Gamma(\mathcal{A},\mathcal{O}_{\mathcal{A}})$ between global sections is exactly the map $\nu$ given in the last proposition.

       As a consequence of the inclusions of Proposition \ref{borderedinclusions}, for surfaces for which the equality $A(\Sigma)=U(\Sigma)$ holds, the three algebras, $A(\Sigma)$, $Sk^{ta}(\Sigma)$, and $U(\Sigma)$ coincide. 

       \begin{cor}
       For any triangulable surface $\Sigma$ with nonempty boundary, the equality $A(\Sigma)=Sk^{ta}(\Sigma)=U(\Sigma)$ holds. 
       \label{borderedAskeinU}
       \end{cor}

       \begin{proof}
       When $\Sigma$ has at least 2 marked points on the boundary or at least 1 puncture together with 1 marked point on the boundary, then $A(\Sigma)=U(\Sigma)$ follows from local acyclicity of the cluster algebra, see Proposition \ref{surfacelocallyacyclic} and Theorem \ref{locallyacyclicau}. When $\Sigma$ has no puncture and 1 marked point on the boundary, let $a$ be its unique invertible boundary variable. Replacing each mutable arc variable $x_\beta$ by $a z_\beta$ divides every homogeneous Ptolemy relation by $a^2$ and gives the coefficient-free relations. The cluster and upper cluster algebras are therefore the corresponding coefficient-free algebras tensored with $\mathbb C[a^{\pm1}]$. Proposition \ref{1markedpointau} gives their equality.\par 
       Now the statement follows from the inclusion $A(\Sigma)\subseteq Sk^{ta}(\Sigma)\subseteq U(\Sigma)$.
       \end{proof}

       \begin{rem}
       The above corollary leaves out only the case that $\Sigma$ has only punctures but no marked points (equivalently, $\Sigma=\Sigma_{g,p}$). For positive genus, $A(\Sigma)\subsetneq U(\Sigma)$ \cite{moon24consequences}. We shall show that $Sk^{ta}(\Sigma)=U(\Sigma)$ for triangulable $\Sigma=\Sigma_{g,p},p\geq 2$ in an alternative manner.
       \end{rem}

		

		\subsection{Canonical bases of cluster algebras}
        The Fock--Goncharov duality conjecture can be generalized to the setting of cluster ensembles\cite{fock09cluster}, where the tropical points $V^{\vee}(\mathbb{Z}^t)$ on the dual $V^{\vee}$ of a cluster variety $V$ are conjectured to parametrize a basis of $\Gamma(V,\mathcal{O}_{V})$. The basis is predicted to be preserved as a set by the cluster modular group of $V$, and to have nonnegative integral structure constants under multiplication. This conjecture reduces to the Fock--Goncharov duality conjecture mentioned in the last subsection in the case of $V=\mathcal{A}_{SL_2,\Sigma}$ and $V^{\vee}=\mathcal{X}_{PGL_2,\Sigma}$.\\
        A significant step towards the general duality conjecture is made by Gross, Hacking, Keel and Kontsevich \cite{gross18canonical} who introduced the theta basis by borrowing techniques from mirror symmetry. In loc. cit., $V^{\vee}(\mathbb{Z}^t)$ parametrizes a basis of a vector space $\mathrm{can}(V)$, with candidate multiplication defined by counts of broken lines. Its subspace $\mathrm{mid}(V)$, with basis indexed by $\Theta\subseteq V^{\vee}(\mathbb{Z}^t)$, is an algebra and admits a $\mathbb{C}$-algebra homomorphism to $\mathrm{up}(V):=\Gamma(V,\mathcal O_V)$. For $V=\mathcal A$, the images of its theta basis contain the set $\bigcup_{\mathbf{s}}\{\prod\limits_{i=1}^n A_i^{m_i}|m_i\geq0,\forall i\in I-I_0\}$ of all the cluster monomials (also called global monomials). For $V$ of type $\mathcal A$, $\mathcal X$, or $\mathcal A_{prin}$, the relevant results of \cite[Theorem 0.3]{gross18canonical} can be summarized as follows:

        \begin{thm}
         (1) Let $\mathrm{can}(V)=\bigoplus\limits_{q\in V^{\vee}(\mathbb{Z}^t)}\mathbb{C}\theta_q$. There is a map $\alpha:V^{\vee}(\mathbb{Z}^t)\times V^{\vee}(\mathbb{Z}^t)\times V^{\vee}(\mathbb{Z}^t)\to \mathbb{Z}_{\geq 0}\cup\{\infty\}$, such that these structure constants define a $\mathbb C$-algebra structure on $\mathrm{can}(V)$ whenever they are finite and, for each $p,q$, only finitely many $r$ satisfy $\alpha(p,q,r)\ne0$.
         
         (2) There is a subset $\Theta\subseteq V^{\vee}(\mathbb Z^t)$ that includes the subset $\Delta^{+}$ of $g$-vectors of global monomials, such that the restriction of $\alpha$ gives their span an associative commutative $\mathbb C$-algebra structure with finite products and nonnegative integral structure constants. This algebra is denoted by the middle cluster algebra $\mathrm{mid}(V)=\bigoplus\limits_{q\in \Theta}\mathbb{C}\theta_q$. 
         
         (3) There is a $\mathbb{C}$-algebra homomorphism $\nu':\mathrm{mid}(V)\to \Gamma(V,\mathcal{O}_V)$, such that, for any $q=\mathbf{g}(m)\in \Delta^{+}$, the theta function $\theta_q$ is mapped to the corresponding global monomial $z^m=\nu'(\theta_q)$.

         (4) The image $\nu'(\theta_q)$ for all $q\in \Theta$ is a universally positive Laurent polynomial, namely a Laurent polynomial with nonnegative integral coefficients in each cluster chart.
        \label{ghkkthm}
        \end{thm}

        In this context, the Fock--Goncharov duality conjecture is stated as follows, as what is called the full Fock--Goncharov conjecture for the cluster variety $V$.
        \begin{conj}[Full Fock--Goncharov conjecture]
        (1) The homomorphism $\nu'$ is injective.
        
        (2) The set $\Theta$ exhausts the tropical points $V^{\vee}(\mathbb{Z}^t)$ on the dual cluster variety, namely $\mathrm{mid}(V)=\mathrm{can}(V)$. 
        
        (3) The homomorphism $\nu'$ is surjective, i.e. $\nu'(\mathrm{can}(V))=\mathrm{up}(V)$.
        \label{fullfgconj}
        \end{conj}

        It has been pointed out by Gross--Hacking--Keel \cite{gross15birational} that the conjecture does not hold in general, while Gross--Hacking--Keel--Kontsevich\cite{gross18canonical} propose several conditions under which the conjecture holds. In the case that $\mathcal{A}$ is associated to a surface, their method, combined with subsequent works, yields important properties of the middle cluster algebra $\mathrm{mid}(\mathcal{A})$. Remarkably, Mandel and Qin\cite{mandel23bracelet} showed that the theta basis, in fact, coincides with the bracelet basis (with few exceptions), which relates the skein algebra to the middle cluster algebra. We summarize notable works concerning the canonical bases of cluster algebras, especially those associated to a surface, in the following theorem.

        \begin{thm}
(1) (\cite[Proposition 8.25]{gross18canonical}, \cite[Proposition 5.23]{mandel23bracelet}) If $\Sigma$ is not a once-punctured closed surface, its cluster complex is not contained in a half-space. Hence $\mathcal A_{prin}$ and its Fock--Goncharov dual have enough global monomials, and the full Fock--Goncharov conjecture holds for $\mathcal A_{prin}$.

(2) (\cite[Theorem 5.24]{mandel23bracelet}) For every triangulable surface, $\nu':\mathrm{mid}(\mathcal A)\to U(\Sigma)$ is injective and $\mathrm{mid}(\mathcal A)=\mathrm{can}(\mathcal A)$.

(3) (\cite[Theorem 8.32(1)]{gross18canonical}) If $\Sigma$ is not a once-punctured closed surface, $\mathrm{can}(\mathcal A)$, and therefore $\mathrm{mid}(\mathcal A)$, is a finitely generated Gorenstein $\mathbb C$-algebra.

(4) (\cite[Theorems 9.4 and 9.6(i)]{mandel23bracelet}) If $\Sigma\ne\Sigma_{1,1}$, the tagged bracelet basis agrees with the theta basis under the maps $\nu$ and $\nu'$. In particular, $Sk^{ta}(\Sigma)=\nu'(\mathrm{mid}(\mathcal A))$ as subalgebras of $U(\Sigma)$.

(5) (\cite[Corollary 5.27]{mandel23bracelet}) If $\Sigma$ is not a once-punctured closed surface, the theta functions form the atomic basis of $\mathrm{mid}(\mathcal A)$ with respect to the cluster atlas. For $\Sigma_{g,1}$ with $g\geq2$, the corresponding statement uses the full tagged triangulation atlas, equivalently the scattering atlas.
\label{ghkkcor}
\end{thm}

The Gorenstein assertion in (3) is applied to the canonical algebra using the enough-global-monomials condition on $\mathcal A_{prin}^{\vee}$ from (1). Parts (2) and (4) identify its image with the tagged skein algebra. Thus this assertion does not depend on the surjectivity of $\nu'$ proved below.
        
        \begin{rem}
         When $\mathcal{A}=\mathcal{A}_{SL_2,\Sigma_{g,p}}$ ($p\geq2$), the conjecture is already known for $\mathcal{A}_{prin}$ and $\mathcal{X}=\mathcal{X}_{PGL_2,\Sigma}$ as a result of the fact that the cluster complex is large. For $\Sigma$ a triangulable surface which is not a $\Sigma_{g,p}$, the full Fock--Goncharov conjecture for $\mathcal{A}$-cluster variety when $G=SL_2$ is implied by the equality of cluster and upper cluster algebra and the results mentioned above. The remaining case addressed here is $\Sigma=\Sigma_{g,p}$ with $p\geq2$; once-punctured closed surfaces are discussed separately in Section 7.
        \end{rem}
        
        In summary, to prove the Fock--Goncharov duality conjecture or the full Fock--Goncharov conjecture for $\Sigma$ not a once-punctured closed surface and $\mathcal{A}=\mathcal{A}_{SL_2,\Sigma}$, it suffices to show the surjectivity of the map $\nu:Sk^{ta}(\Sigma)\to U(\Sigma)$, which is equivalent to the surjectivity of $\nu': \mathrm{mid}(\mathcal A)\to \Gamma(\mathcal{A},\mathcal{O}_{\mathcal{A}})$.

        \begin{prop}
        If $\Sigma$ is not a once-punctured closed surface and $\mathcal{A}=\mathcal{A}_{SL_2,\Sigma}$, the Fock--Goncharov duality conjecture (in the sense of Conjecture \ref{fgdualityconj}) and the full Fock--Goncharov conjecture (in the sense of Conjecture \ref{fullfgconj}) are equivalent to the surjectivity of the map $\nu$, or the equality $Sk^{ta}(\Sigma)=U(\Sigma)$. (By abuse of notation, we identify the skein algebras with their image in the upper cluster algebra.)
        \label{conjs}
        \end{prop}

		\section{Birational geometry of skein and cluster algebras}

Throughout this section, let $\Sigma$ be a triangulable marked surface. Instead of directly comparing skein and cluster algebras, it is useful to compare their spectra with the cluster variety. We first distinguish a common open subscheme from the whole spectrum of the Roger--Yang algebra.

\subsection{A common open subscheme of the two spectra}
Let
\[
 \begin{gathered}
 R=Sk^{ta}(\Sigma),\qquad B=Sk^{RY}(\Sigma)=R[v_1^{\pm1},\ldots,v_p^{\pm1}],\\
 X=\operatorname{Spec}R,\qquad Y=\operatorname{Spec}B.
 \end{gathered}
\]
The two algebras are finitely generated domains in the same rational function field. Write $\phi:Y\to X$ for the morphism induced by the inclusion. Let $\mathcal S\subseteq X$ be the largest open subscheme on which all the rational functions $v_i^{\pm1}$ are regular.

\begin{prop}\label{RYtaggedbirational}\label{commonlocus}
There is a canonical open immersion $\mathcal S\to Y$. The morphism $\phi$ induces an isomorphism $\phi^{-1}(\mathcal S)\cong\mathcal S$, and its restriction is the inclusion $\mathcal S\to X$.
\end{prop}
\begin{proof}
For each point of $\mathcal S$, choose a principal neighborhood $D_X(f)\subseteq\mathcal S$. Since every $v_i^{\pm1}$ belongs to $R_f$, we have $B_f=R_f$. Thus $\phi^{-1}(D_X(f))=D_Y(f)\cong D_X(f)$. These isomorphisms agree in the common rational function field and therefore glue.
\end{proof}

This does not assert that $Y\to X$ is an open immersion: points of $Y$ outside $\mathcal S$ may remain. In particular, adjoining rational functions is not in general the same operation as restricting to their common regular locus.

\subsection{Cluster varieties and extension of regular functions}
We use the cluster variety obtained by gluing seed tori on their maximal common open subsets in the rational function field. For a seed with cluster variables $x_1,\ldots,x_n$, put $f=x_1\cdots x_n$, including the invertible frozen variables when present. The inclusions
\[
 \mathbb C[x_1,\ldots,x_n]\subseteq R\subseteq
 \mathbb C[x_1^{\pm1},\ldots,x_n^{\pm1}]
\]
give $R_f=\mathbb C[x_1^{\pm1},\ldots,x_n^{\pm1}]$. Hence the seed torus is $D_X(f)$, and the seed inclusions glue to an open immersion $\mathcal A\to X$. In particular, this version of $\mathcal A$ is integral, separated and Noetherian. Its ring of global sections is the upper cluster algebra. The same localization observation applies to intermediate cluster algebras in \cite{beyer25deep}. The extension criterion below is the commutative-algebra principle underlying \cite[Lemmas 1.2.1 and 4.4.1]{matherne15upper}.

\begin{lem}\label{basiclemma}
Let $R$ be a Noetherian domain and $V$ a nonempty open subscheme of $X=\operatorname{Spec}R$.
If $R$ satisfies Serre's condition $S_2$ and $X\setminus V$ has codimension at least two, then $R=\Gamma(V,\mathcal O_V)$.
Conversely, if $R=\Gamma(V,\mathcal O_V)$, then $X\setminus V$ contains no point of height one.
In particular, for a normal Noetherian domain these two conditions are equivalent.
\end{lem}
\begin{proof}
For an $S_2$ domain, the intersection of $R_{\mathfrak p}$ over the height-one primes $\mathfrak p$, taken in $\operatorname{Frac}R$, equals $R$. A regular function on $V$ belongs to all these local rings when $V$ contains all height-one points, proving the first assertion.

For the converse, suppose that a height-one prime $\mathfrak p$ is absent from $V$. Choose $0\ne a\in\mathfrak p$. The prime $\mathfrak p$ is minimal over $(a)$ and hence is an associated prime of $R/(a)$. There is therefore $b\notin(a)$ with $(a):b=\mathfrak p$. The rational function $b/a$ does not belong to $R$, but is regular on $X\setminus V(\mathfrak p)$: if $\mathfrak q$ does not contain $\mathfrak p$, choose $c\in\mathfrak p\setminus\mathfrak q$ and write $cb=ad$, so $b/a=d/c\in R_{\mathfrak q}$. Since the closed set $X\setminus V$ contains $V(\mathfrak p)$, this function is regular on $V$, a contradiction.
\end{proof}

\begin{prop}\label{borderedgeometry}
For a marked surface with boundary for which $A(\Sigma)=Sk^{ta}(\Sigma)=U(\Sigma)$, the open immersion $\mathcal A\to\operatorname{Spec}(Sk^{ta}(\Sigma))$ has complement of codimension at least two. This applies to the locally acyclic cases of Proposition~\ref{surfacelocallyacyclic}, and to the coefficient-free once-marked case of Proposition~\ref{1markedpointau}.
\end{prop}
\begin{proof}
The upper cluster algebra is an intersection of Laurent polynomial rings in a common field and hence is normal. The assertion follows from Lemma~\ref{basiclemma} and the assumed equality of global sections.
\end{proof}

\subsection{Normalization and the localized cluster variety}
The following comparison uses the dimension estimate of Section~4 and the main theorem of Section~5. Neither of those proofs uses the present proposition.

\begin{prop}\label{normalizationcomparison}
Let $\Sigma=\Sigma_{g,p}$ be triangulable with $p\geq2$. Let $\mathring{\mathcal A}=\mathcal A\cap\mathcal S$ and let $\overline B$ be the integral closure of $B$ in its fraction field. Then $\mathring{\mathcal A}\to Y$ is an open immersion with complement of codimension at least two, and
\[
 \Gamma(\mathring{\mathcal A},\mathcal O_{\mathring{\mathcal A}})
 =\Gamma(\mathcal S,\mathcal O_{\mathcal S})=\overline B.
\]
Consequently, equality with $B$ is equivalent to normality of $B$.
\end{prop}
\begin{proof}
Put $d=6g-6+3p=\dim R=\dim B$. Theorem~\ref{mainthmproof} shows that $X$ is normal and $X\setminus\mathcal A$ has codimension at least two. Proposition~\ref{commonlocus} already gives the required open immersions into $Y$.

For a nonempty set $J$ of punctures, let $I_J\subseteq R$ be generated by the tagged arcs incident to at least one puncture of $J$. Proposition~\ref{multipuncturedimension} gives
\begin{equation}\label{stratumdimensionbound}
 \dim(R/I_J)+|J|\leq d-2.
\end{equation}
Stratify $Y$ according to the set $J$ of punctures at which every incident tagged arc vanishes. At a point of this stratum and for each $i\notin J$, choose a plain arc $a_i$ incident to the $i$th puncture which does not vanish. This is possible because different taggings differ by units in $B$. If $m_i\in\{1,2\}$ is its number of endpoints at that puncture, the fully changed arc $b_i=v_i^{m_i}a_i$ also does not vanish. The generalized tagged curves $v_i a_i$ and $v_i^{m_i-1}a_i$ belong to $R$ by Definition--Proposition~\ref{taggedskein}, so on $D_X(a_i b_i)$ we have
\[
 v_i=\frac{v_i a_i}{a_i},\qquad
 v_i^{-1}=\frac{v_i^{m_i-1}a_i}{b_i}.
\]
With $f=\prod_{i\notin J}a_i b_i$, this gives
\[
 B_f=R_f[v_j^{\pm1}:j\in J].
\]
The part of the stratum in $D_Y(f)$ lies in $V_Y(I_J)$, whose coordinate ring there is a quotient of
$(R/I_J)_f[t_j^{\pm1}:j\in J]$. Its dimension is at most $d-2$ by~\eqref{stratumdimensionbound}. Each stratum is locally closed in a Noetherian scheme, so a finite number of these opens suffices. The stratum with $J=\varnothing$ lies in $\mathcal S$. There are only finitely many sets $J$, and therefore $Y\setminus\mathcal S$ has dimension at most $d-2$. Together with the bound for $X\setminus\mathcal A$, this proves the codimension assertion for $Y\setminus\mathring{\mathcal A}$.

The normalization $\overline Y=\operatorname{Spec}\overline B\to Y$ is finite because $B$ is a finitely generated algebra over $\mathbb C$. It is an isomorphism over the normal open subscheme $\mathcal S$, and the inverse image of its complement has codimension at least two. The same holds for $\mathring{\mathcal A}$, which is smooth. Applying Lemma~\ref{basiclemma} to $\overline Y$ proves both equalities.
\end{proof}

\section{Estimation of Krull dimension}

\subsection{Geometric intersection numbers and skein products}
Let $CS(\Sigma)$ denote the set of curve systems on $\Sigma$, and let
$TCS(\Sigma)$ denote the set of tagged curve systems. Tags are forgotten
when taking geometric intersection numbers. If $K$ is a simple closed
curve, write $I(\mu,K)$ for the minimum number of transverse intersections
between representatives of $\mu$ and $K$. For a puncture $Q$, let
$e_Q(\mu)$ denote the number of ends of $\mu$ at $Q$.

\begin{prop}
\label{intersectionresolution}
Let $\gamma$ be a curve system occurring after resolving a product of
(tagged) curve systems $\alpha\beta$ by the skein relations. For every
simple closed curve $K$ and every puncture $Q$,
\[
 I(\gamma,K)\leq I(\alpha,K)+I(\beta,K),\qquad
 e_Q(\gamma)\leq e_Q(\alpha)+e_Q(\beta).
\]
\end{prop}
\begin{proof}
Fix a representative of $K$. Choose representatives of each factor
realizing its minimum intersection with $K$, and perturb them so that no
crossing to be resolved lies on $K$. Resolve crossings in disks disjoint
from $K$, and resolve at punctures in sufficiently small neighborhoods
also disjoint from $K$. These operations create no new intersections with
$K$ and no new ends at a puncture. Removing an inessential component does
not increase either quantity. Taking minimum intersection numbers gives
the assertion. This argument is applied separately for each test curve;
no simultaneous tautness assertion is needed.
\end{proof}

\subsection{Dehn--Thurston coordinates and polynomial growth}
Suppose $2g-2+q>0$. Replace the punctures of $\Sigma_{g,q}$ by boundary
circles and choose a pants decomposition with $r=3g-3+q$ interior cuffs.
The Dehn--Thurston coordinates of a curve system consist of $r$
intersection coordinates, $r$ twist coordinates, and $q$ endpoint
coordinates. Equivalently, they lie in
\[
 (\mathbb Z^2/\{\pm1\})^r\times\mathbb Z_{\geq0}^{q}.
\]
We choose the representative $(n_i,t_i)$ of each pair with $n_i\geq0$;
when $n_i=0$, the twist coordinate records the nonnegative multiplicity
of the corresponding cuff. In particular, a bounded coordinate box
contains $O(N^{6g-6+3q})$ curve systems. Parity conditions only reduce this
number. Put
\[
 \|\Pi(\mu)\|=
 \max\bigl(\{n_i(\mu),|t_i(\mu)|:1\leq i\leq r\}
               \cup\{e_Q(\mu):Q\text{ a puncture}\}\bigr),
\]
with the maximum of an empty collection taken to be zero.

\begin{prop}
\label{luodtintersection}
For a fixed pants decomposition, there is a finite collection of simple
closed curves $K_1,\ldots,K_s$ such that every Dehn--Thurston coordinate
pair is a continuous homogeneous function of
\[
 \bigl(I(\mu,K_1),\ldots,I(\mu,K_s),
       (e_Q(\mu))_{Q}\bigr).
\]
\end{prop}
\begin{proof}
The local reconstruction formulas on a four-holed sphere and a one-holed
torus in \cite[Proposition 4.4]{luo04dt} recover the coordinate pair at
each interior cuff from intersections with finitely many curves and the
boundary weights. Apply these formulas to all cuffs. The weights at the
original boundary circles are exactly the endpoint counts $e_Q$.
\end{proof}

\begin{lem}
\label{dtintersection}
There are $K_1,\ldots,K_s$ and $C>0$, depending only on the surface and
coordinate choices, such that
\[
 \|\Pi(\mu)\|\leq
 C\left(\sum_{j=1}^s I(\mu,K_j)+\sum_Q e_Q(\mu)\right).
\]
Consequently, the number of tagged curve systems for which the expression
in parentheses is at most $N$ is $O((N+1)^{6g-6+3q})$.
\end{lem}
\begin{proof}
For each continuous homogeneous coordinate recovery function in
Proposition~\ref{luodtintersection}, its norm is bounded on the compact
unit simplex in the nonnegative input coordinates. Homogeneity then gives
a linear bound on all inputs. Taking the maximum of the finitely many
bounds proves the first assertion; compare the same argument in
\cite[Lemma 9.6]{bloomquist23degenerations}. Integral coordinate counting
proves the second assertion. There are at most $2^q$ choices of tags for
an underlying curve system, so tagging does not change the exponent.
\end{proof}

\subsection{Filling punctures and dimension bounds}
In this subsection, let $\Sigma=\Sigma_{g,p}$, where $g=0,p\geq4$ or
$g\geq1,p\geq1$, and put
\[
 R=Sk^{ta}(\Sigma),\qquad d=6g-6+3p.
\]
The finite generation and domain statements in
Theorem~\ref{skeinalgebraic}(1), together with the common fraction field
in Proposition~\ref{inclusions}, give $\dim R=d$. This dimension
calculation does not use the equality $R=U(\Sigma)$.

\begin{lem}
\label{fillingspanning}
Let $J$ be a nonempty set of punctures, and let $I_J\subseteq R$ be the
ideal generated by all tagged arcs with an end in $J$. Set
$B_J=R/\sqrt{I_J}$ and let $\Sigma_J=\Sigma_{g,p-|J|}$ be obtained by
filling the punctures in $J$. The following assertions hold.
\begin{enumerate}
\item Every generalized tagged arc with an end in $J$ has zero image
in $B_J$.
\item Moving a strand across a point of $J$ changes its image in $B_J$
by a sign.
\item The vector space $B_J$ is spanned by the images of at most one
chosen lift of each reduced tagged curve system on $\Sigma_J$.
\end{enumerate}
Here reduced means that no component is contractible and no loop is
peripheral, as in the definition of a curve system. The last assertion
is a spanning statement; it does not assert linear independence or an
identification with a skein algebra of $\Sigma_J$.
\end{lem}
\begin{proof}
Resolve the self-intersections of a generalized arc and apply the digon
relations. An admissible tagged arc ending in $J$ is a generator of
$I_J$. If a simple returning arc has different tags at its two ends,
its square is the product of the same underlying arc with both ends
plain and with both ends notched, and hence lies in $I_J$. An arc
bounding a once-punctured monogon is expressed by the corresponding
product of tagged radii; if its base lies in $J$, both radii end in $J$.
Contractible arcs vanish. These reductions prove (1).

For (2), the digon relation expresses the sum of the two diagrams before
and after a strand is moved across a puncture $P\in J$ as a diagram
with ends at $P$. Its image is zero by (1). Thus the two images are
negatives of one another.

Expand an element of $R$ in tagged curve systems and discard those with
an end in $J$. Regard each remaining system as a diagram on $\Sigma_J$.
If a component becomes contractible or a loop becomes peripheral, move
it into a disk or punctured disk in $\Sigma_J$ and apply the framing or
contractible-arc relation. The isotopy can be made generic with respect
to the filled points: away from finitely many strand crossings of these
points it is an isotopy on $\Sigma$. Each such crossing changes only a
sign by (2). Thus the reduction on the filled surface is valid in $B_J$
up to a scalar. Finally, two lifts of the same reduced tagged curve
system are joined by such a generic isotopy and have proportional
images. This proves (3). Possible additional relations among the lifts
can only decrease the dimension of their span.
\end{proof}

\begin{prop}
\label{multipuncturedimension}
With the preceding notation, set $q=p-|J|$. Then
\[
 \dim(R/I_J)=\dim B_J\leq D(g,q),
\]
where
\[
 D(g,q)=
 \begin{cases}
 6g-6+3q,&2g-2+q>0,\\
 2,&(g,q)=(1,0),\\
 1,&(g,q)=(0,2),\\
 0,&(g,q)=(0,0)\text{ or }(0,1).
 \end{cases}
\]
In particular, when $p\geq2$ in the range above,
$\operatorname{ht} I_J\geq |J|+2$.
\end{prop}
\begin{proof}
Choose a finite-dimensional generating subspace $V\subseteq R$
containing $1$ and spanned by finitely many tagged curve systems. Let
$F_N$ be the image of $V^N$ in $B_J$. Assume first that
$2g-2+q>0$, and choose the test curves of Lemma~\ref{dtintersection}
on $\Sigma_J$. Represent these curves away from the filled points and
regard them also as curves on $\Sigma$. Proposition~\ref{intersectionresolution}
gives a constant $c$, independent of $N$, such that every system arising
from resolution of a word of length at most $N$ has total test
intersection and total endpoint count at the remaining punctures at
most $cN$.

Filling punctures cannot increase geometric intersection numbers with
the chosen test curves. Deleting inessential components cannot increase
these numbers or endpoint counts either. By Lemma~\ref{fillingspanning},
$F_N$ is therefore spanned by chosen lifts of reduced tagged curve systems
on $\Sigma_J$ satisfying the same bound. Lemma~\ref{dtintersection}
implies
\[
 \dim_{\mathbb C}F_N=O((N+1)^{6g-6+3q}).
\]
Since Gelfand--Kirillov dimension equals Krull dimension for a finitely
generated commutative algebra, this proves the asserted bound in the
stable cases; see \cite{mcconnellnoncommutative}.

For a closed torus, a reduced nonempty curve system is a collection of
parallel essential loops. Its slope and multiplicity are encoded by
an integral pair $(u,v)$ up to simultaneous sign. Intersections with a
meridian and a longitude bound $|u|$ and $|v|$, giving $O((N+1)^2)$
possibilities. On a sphere with two punctures, every reduced nonempty
system consists of parallel copies of the unique arc joining the two
punctures; its multiplicity is bounded by its endpoint count, giving
$O(N+1)$ possibilities. A sphere with at most one puncture has only the
empty reduced system. The same spanning argument proves the remaining
bounds.

For an affine domain over $\mathbb C$, the dimension formula gives
$\operatorname{ht} I_J=d-\dim(R/I_J)$ whenever $I_J$ is proper. The
assertion is automatic if $I_J=R$. Direct substitution in the displayed
values of $D(g,q)$ gives $d-D(g,q)\geq |J|+2$ for $p\geq2$.
\end{proof}

\begin{lem}
\label{dimensionlemma}
Let $P$ be a puncture and let $I_P$ be the ideal of $R$ generated by all
tagged arcs incident to $P$. Then $\operatorname{ht} I_P\geq3$.
Consequently, the ideal $I_0$ generated by all tagged arcs also has
height at least three.
\end{lem}
\begin{proof}
Except when $(g,p)=(1,1)$, filling $P$ gives a surface with
$2g-2+(p-1)>0$. Proposition~\ref{multipuncturedimension} therefore gives
$\dim R/I_P\leq d-3$.

For the once-punctured torus, let $\ell$ be an essential simple loop.
Choose a simple returning arc $\alpha$ based at $P$ whose two puncture
smoothings are isotopic to $\ell$. The puncture skein relation gives
$\alpha_{+-}=2\ell$, where $\alpha_{+-}$ has one plain and one notched
end. Hence
\[
 4\ell^2=\alpha_{+-}^2=\alpha_{++}\alpha_{--}\in I_P.
\]
Thus every essential simple loop vanishes in $R/\sqrt{I_P}$. The tagged
arcs also vanish, so this quotient is spanned by $1$ and has dimension
at most zero. Since $\dim R=3$, the height bound follows. Finally,
$I_P\subseteq I_0$ gives the last assertion.
\end{proof}

\section{Proof of the main theorem \texorpdfstring{$Sk^{ta}(\Sigma)=U(\Sigma)$}{Skta = U}}

The proof first establishes equality after inverting one tagged arc and
then uses the height estimate to extend regular functions. In particular,
it does not require a global open immersion
$\operatorname{Spec}Sk^{RY}(\Sigma)\to\operatorname{Spec}Sk^{ta}(\Sigma)$.

\begin{lem}
\label{onearcfreezing}
Let $\Sigma=\Sigma_{g,p}$, where $g=0,p\geq4$ or $g\geq1,p\geq1$.
Choose a tagged arc $\alpha$ in a seed, write $a=x_\alpha$, and denote
by $A^{[\alpha]}$ and $U^{[\alpha]}$ the cluster and upper cluster
algebras obtained by freezing $a$, with $a$ inverted. Then
\[
 A^{[\alpha]}=U^{[\alpha]}.
\]
For $p=1$, the seed is taken in the tag component containing $\alpha$.
\end{lem}
\begin{proof}
Change tags at punctures, if necessary, so that $\alpha$ is plain, and
cut along its underlying arc. By the cutting--freezing correspondence,
the mutable quiver of the frozen seed is that of the cut surface.

If the cut surface is connected, it has two boundary marked points.
The covering-pair induction for such surfaces in
\cite[Theorem 10.6]{muller13acyclic} applies to this mutable quiver with
its frozen coefficients. It is locally acyclic, so
$A^{[\alpha]}=U^{[\alpha]}$ by Theorem~\ref{locallyacyclicau}.

Suppose that cutting disconnects the surface into $\Sigma_1$ and
$\Sigma_2$. Each component has one boundary component and one boundary
marked point. Here $\alpha$ is a separating arc with coincident endpoints,
so it cannot be a folded radius. The two component seeds therefore have
their boundary variables identified with $a$, and mutations in either
component preserve this identification. A punctured component has a Louise mutable quiver, by the
proof of \cite[Proposition 12]{canakci15cluster}. The Louise induction
is a condition on the mutable quiver and remains valid with the frozen
coefficient $a$; see \cite[Definition 2.3 and Proposition 2.6]{mullerspeyer16louise}. If both components are punctured, their disjoint union
is again Louise, so the result follows from local acyclicity.

It remains to handle any unpunctured components. Such a component has
positive genus, since $\alpha$ is not contractible. It has a single
boundary segment with variable $a$. In every seed on this component,
put $z_\beta=a^{-1}x_\beta$ for each mutable arc $\beta$. Each surface
exchange relation is homogeneous of degree two when all arc and
boundary-segment variables have degree one. Dividing it by $a^2$
therefore gives exactly the coefficient-free exchange relation in the
$z$ variables. Theorem~5 of \cite{canakci15cluster} applies to this
coefficient-free cluster algebra and gives $A_i^0=U_i^0$.

These normalized variables on different components are algebraically
independent. If one component is unpunctured and the other is punctured,
write $A_2,U_2$ for the latter component with coefficient $a$; then
\[
 A^{[\alpha]}=A_1^0\otimes_{\mathbb C}A_2,
 \qquad
 U^{[\alpha]}=U_1^0\otimes_{\mathbb C}U_2.
\]
If both are unpunctured, the corresponding equalities are
\[
 A^{[\alpha]}=A_1^0\otimes_{\mathbb C}A_2^0
                    \otimes_{\mathbb C}\mathbb C[a^{\pm1}],
 \qquad
 U^{[\alpha]}=U_1^0\otimes_{\mathbb C}U_2^0
                    \otimes_{\mathbb C}\mathbb C[a^{\pm1}].
\]
For the upper-algebra equalities, intersect the Laurent rings over all
pairs of component seeds. Intersections of subspaces commute with tensor
products over a field: for an element under consideration, expand in a
finite linearly independent set in one factor and compare coefficients.
Apply this observation successively to the two independent sets of
variables. Thus each upper algebra is the indicated tensor product.
The equalities for the component algebras prove the lemma.
\end{proof}

\begin{thm}
\label{mainthmproof}
Let $\Sigma=\Sigma_{g,p}$ with $g=0,p\geq4$ or $g\geq1,p\geq2$.
Then the natural inclusion is an equality
\[
 Sk^{ta}(\Sigma)=U(\Sigma).
\]
Moreover, the open immersion
$\mathcal A\to\operatorname{Spec}Sk^{ta}(\Sigma)$ has complement of
codimension at least two.
\end{thm}
\begin{proof}
Put $R=Sk^{ta}(\Sigma)$ and $X=\operatorname{Spec}R$. We use finite
generation and the $S_2$ property from Theorem~\ref{skeinalgebraic};
these are independent inputs established before the present theorem.
For every tagged arc $\alpha$, its variable $a=x_\alpha$ is a cluster
variable because $p\geq2$. The freezing inclusions and
$A(\Sigma)\subseteq R\subseteq U(\Sigma)$ give
\[
 A^{[\alpha]}\subseteq A(\Sigma)[a^{-1}]
 \subseteq R[a^{-1}]\subseteq U(\Sigma)[a^{-1}]
 \subseteq U^{[\alpha]}.
\]
Lemma~\ref{onearcfreezing} makes all these inclusions equalities.
In particular, every $f\in U(\Sigma)$ is regular on $D_R(a)$.

Let
\[
 X^\circ=\bigcup_{\alpha}D_R(x_\alpha),
\]
where the union ranges over all tagged arcs. Its complement is $V(I_0)$,
which has codimension at least three by Lemma~\ref{dimensionlemma}.
The $S_2$ extension property for the Noetherian domain $R$ gives
\[
 R=\Gamma(X^\circ,\mathcal O_{X^\circ}).
\]
Since every element of $U(\Sigma)$ is regular on $X^\circ$, this proves
$U(\Sigma)\subseteq R$, and the opposite inclusion is already known.

Finally, $U(\Sigma)$ is an intersection of Laurent polynomial rings in
the common fraction field and is therefore normal. Thus $R$ is normal.
For each seed, localizing $R$ at the product of its cluster variables
gives precisely its Laurent polynomial ring. Consequently the seed tori
are principal open subsets of $X$, and, with the gluing convention of
Subsection~2.4, their union is the given open subscheme $\mathcal A$.
We have just proved
\[
 \Gamma(\mathcal A,\mathcal O_{\mathcal A})=U(\Sigma)=R.
\]
Lemma~\ref{basiclemma} now implies that $X\setminus\mathcal A$ has
codimension at least two.
\end{proof}

\begin{rem}
The filling argument in Section~4 supplies spanning sets and dimension
upper bounds. It does not identify a quotient by arc classes with the
skein algebra of a filled surface. In particular, no assertion about
non-affineness or the nonemptiness of $X\setminus\mathcal A$ follows
from that argument alone.
\end{rem}

\section{Applications of the main theorem}

The comparison with the cluster variety gives algebraic properties of the
tagged skein algebra.  For the Roger--Yang algebra, Proposition
\ref{normalizationcomparison} identifies the ring of functions on the
localized cluster variety with its normalization.  To replace normalization
by equality, we give an independent normality argument.  The four-punctured
sphere is an exception, as the explicit example below shows.

\subsection{Normality of the Roger--Yang skein algebra}

Write $B_{g,p}=Sk^{RY}(\Sigma_{g,p})$ and
\[
 S_{g,p}=B_{g,p}/(v_1-1,\ldots,v_p-1).
\]
The specialization is at $q^{1/2}=1$ throughout this section.  We use the
normalization of the puncture relation in which its coefficient is $v_j$;
thus the coefficient in $S_{g,p}$ is $1$.

\begin{thm}\label{RYnormality}
The classical Roger--Yang skein algebra $B_{g,p}$ is a finitely generated
normal domain if $g\geq1$ and $p\geq1$, or if $g=0$ and $p\geq5$.
\end{thm}

We first construct the affine model used in the proof.  Put $V=\mathbb C^2$,
$G=SL_2(\mathbb C)$, and
\[
 \omega(u,v)=u_1v_2-u_2v_1,\qquad
 N(w)u=\omega(w,u)w,\qquad P(w)=-I+N(w).
\]
Then $N(w)^2=0$, $P(w)\in G$, and $\operatorname{tr}P(w)=-2$.
With $[A,B]=ABA^{-1}B^{-1}$, define the affine scheme
\begin{equation}\label{RYmatrixmodel}
 Z_{g,p}=\left\{(\mathbf A,\mathbf B,\mathbf w)\in G^{2g}\times V^p:
 \prod_{i=1}^g[A_i,B_i]\prod_{j=1}^pP(w_j)=I\right\}.
\end{equation}
Here $G$ acts by simultaneous conjugation on the matrices and by its
standard action on the vectors.  Zero vectors are allowed.

\begin{lem}\label{RYmatrixnormal}
In the range of Theorem \ref{RYnormality}, the scheme $Z_{g,p}$ is normal
and is a local complete intersection of pure dimension $6g+2p-3$.
\end{lem}

\begin{proof}
Let $\mu:G^{2g}\times V^p\to G$ be the product in
\eqref{RYmatrixmodel}, and set $D=6g+2p-3$.  We compute the critical
locus of the fiber $\mu^{-1}(I)$.  Using the trace pairing to identify
$\mathfrak{sl}_2$ with its dual, the annihilator of $d\mu$ at a point
of this fiber is
\begin{equation}\label{RYdifferential}
 \{H\in\mathfrak{sl}_2:HA_i=A_iH,\ HB_i=B_iH\text{ for all }i,
                  \quad Hw_j=0\text{ for all }j\}.
\end{equation}
We include the calculation.  For a puncture factor,
\[
 \operatorname{tr}\bigl(H\,dP_w(\dot w)P(w)^{-1}\bigr)
 =\omega\bigl(2Hw+\omega(w,Hw)w,\dot w\bigr).
\]
Its vanishing for every $\dot w$ implies $Hw=0$: pairing the vector in
parentheses with $w$ first gives $\omega(w,Hw)=0$.  In particular,
$H$ then commutes with $P(w)$.  For a commutator $C=[A,B]$, the
right-translated derivatives in the directions $\dot A=XA$ and
$\dot B=YB$ are respectively
\[
 X-\operatorname{Ad}_{ABA^{-1}}X,\qquad
 \operatorname{Ad}_A Y-\operatorname{Ad}_C Y.
\]
Their annihilator consists precisely of the $H$ commuting with $A$
and $B$.  Indeed, the first equation says that
$K=\operatorname{Ad}_{A^{-1}}H$ commutes with $B$; the second then says
$K=\operatorname{Ad}_{BA}K$, hence that $K$ also commutes with $A$.
Finally, in differentiating the product, the transported annihilator is
unchanged after each factor, because it commutes with that factor.  This
proves \eqref{RYdifferential}.

If a nonzero $H$ in \eqref{RYdifferential} is semisimple, every $w_j$
is zero and all handle matrices belong to one maximal torus.  The
product equation requires $p$ even, and this locus has dimension at
most $2g+2$.  If $H$ is nilpotent, choose a nonzero vector $e$ spanning
$\ker H$.  Every handle matrix lies in the centralizer of $H$, whose
elements have the form $\pm(I+aH)$, and $w_j=c_je$.  These matrices
commute, so the product equation becomes
\[
 (-1)^p\left(I-\Bigl(\sum_j c_j^2\Bigr)N(e)\right)=I.
\]
Thus $p$ is even and $\sum_jc_j^2=0$.  Choosing the line $\mathbb Ce$
adds one parameter; the resulting dimension bound is
$1+2g+(p-1)=2g+p$.  This also bounds the semisimple locus when $p\geq2$.

For odd $p$ the fiber has no critical points.  For even $p$ in the
stated range, the critical locus has dimension at most $2g+p$, and
\[
 D-(2g+p)=4g+p-3\geq3.
\]
Every component of the fiber has dimension at least $D$, since a fiber
over the smooth three-dimensional group $G$ is locally cut out by
three equations.  No component can be contained in the critical locus;
at its general point the fiber is smooth of dimension $D$.  Hence the
fiber has pure dimension $D$ and is a local complete intersection.
It is Cohen--Macaulay and generically reduced, hence reduced.  The
critical-locus bound implies regularity in codimension one.  Serre's
criterion now proves normality.
\end{proof}

\begin{lem}\label{RYinvariantmodel}
In the range of Theorem \ref{RYnormality}, curve evaluation induces an
isomorphism
\[
 S_{g,p}\ \simeq\ \mathbb C[Z_{g,p}]^G.
\]
\end{lem}

\begin{proof}
We specify the curve evaluation, including its signs.  Replace the
punctures by small boundary circles.  For the unit tangent bundle, let
$z$ be a positive fiber loop and let $d_j$ be the tangent lift of the
$j$-th oriented boundary circle.  Handle lifts may be chosen with
\[
 \prod_i[a_i,b_i]\prod_j d_j=z^{2-2g-p},\qquad z\text{ central}.
\]
A point of $Z_{g,p}$ defines a twisted local system by
\[
 a_i\longmapsto A_i,\quad b_i\longmapsto B_i,\quad
 z\longmapsto-I,\quad d_j\longmapsto I-N(w_j)=-P(w_j).
\]
The displayed relation holds because both sides evaluate to
$(-1)^pI$.  Moreover, $d_jw_j=w_j$, so $w_j$ is a monodromy-invariant boundary
vector; the affine model allows this vector to vanish.

For arcs we use the good lift of \cite[Section 8.1, Definition 8]{costantino22stated}:
bend each endpoint tangent to the positive boundary tangent before
taking its tangent lift.  The evaluation of an arc $\alpha:i\to j$
and of an immersed loop $\gamma$ is respectively
\[
 F_\alpha=\omega\bigl(w_j,\rho(\widehat\alpha)w_i\bigr),
 \qquad F_\gamma=\operatorname{tr}\rho(\widetilde\gamma).
\]
The reversal identity
$\rho(\widehat{\alpha^{-1}})=-\rho(\widehat\alpha)^{-1}$
proves orientation independence; see
\cite[Lemma 8.7]{costantino22stated}.  The ordinary skein relations
follow from \cite[Theorem 8.12]{costantino22stated}, extended by
multilinearity from stated basis vectors to arbitrary endpoint vectors.

For completeness, the additional puncture relation is a separate
local calculation.  In a boundary collar its two smoothings have
transport $I$ and $-d_j$: the complementary passage introduces one
full tangent-fiber turn.  If the two remaining transported vectors
are $u$ and $v$, their sum evaluates to
\[
 \omega(v,u)-\omega(v,d_ju)
 =\omega\bigl(v,N(w_j)u\bigr)
 =\omega(v,w_j)\omega(w_j,u),
\]
which is the normalized puncture relation.  A peripheral loop has
value $2$, and a contractible loop has value $\operatorname{tr}(-I)=-2$.
Sliding an endpoint around its puncture preserves the value because
$d_jw_j=w_j$.  Thus evaluation gives an algebra homomorphism
\begin{equation}\label{RYevaluationmap}
 S_{g,p}\longrightarrow\mathbb C[Z_{g,p}]^G.
\end{equation}
All these functions are polynomial in the matrix entries and the
vectors: inversion on $SL_2$ is polynomial.  Consequently the
construction applies on all of $Z_{g,p}$, including zero decorations.

The mixed invariant theorem
\cite[Theorem 6.4.1]{muller13decorated} says that the invariants on
$G^{2g}\times V^p$ are generated by traces of matrix words and the
contractions $\omega(w_j,Ww_k)$.  Since $G$ is reductive over
$\mathbb C$, taking invariants is exact, so their restrictions generate
$\mathbb C[Z_{g,p}]^G$.  Each is, up to its fixed good-lift sign, the
evaluation of an immersed loop or arc.  Hence
\eqref{RYevaluationmap} is surjective.  This use of invariant theory
does not assume a prior identification of the Roger--Yang algebra
with a decorated character algebra.

By \cite[Theorems 5.3 and 9.1]{bloomquist23degenerations}, the algebra
$S_{g,p}$ is a finitely generated domain of dimension $6g-6+2p$ in
the range under consideration.  We show that the target of
\eqref{RYevaluationmap} has at least this dimension.  By Lemma
\ref{RYmatrixnormal}, $Z_{g,p}$ has pure dimension $D=6g+2p-3$.
It has a nonempty open set of closed orbits of dimension three, as
the following points show.

If $g\geq1$ and $p$ is odd, put every $w_j=0$, take
\[
 A_1=\begin{pmatrix}\mathrm i&0\\0&-\mathrm i\end{pmatrix},
 \qquad B_1=\begin{pmatrix}0&1\\-1&0\end{pmatrix},
\]
and set the other handle matrices equal to $I$.  The first commutator
is $-I$, and the first handle pair is irreducible.  If $g\geq1$ and
$p$ is even, take $w_1=e_1$, $w_2=e_2$, and the other vectors zero.
The matrix $K=(P(e_1)P(e_2))^{-1}$ has distinct eigenvalues and trace
$1$.  It is a commutator of an irreducible pair: after diagonalizing
$K=\operatorname{diag}(\lambda,\lambda^{-1})$, choose
$A_1=\operatorname{diag}(s,s^{-1})$ with $s^2=\lambda$ and the same
antidiagonal $B_1$ as above.  Conjugating this pair back gives the
required point.  An irreducible matrix pair has a closed conjugacy
orbit and finite stabilizer, and the same holds on its open
irreducible locus with the decorations included.

For $g=0$ and odd $p\geq5$, take
\[
 w_1=(1,0),\quad w_2=(0,2),\quad w_3=(1,-2),
 \quad w_j=0\ (j>3).
\]
The first three $P$-matrices have product $I$.  For even $p\geq6$,
take $(w_1,w_2,w_3,w_4)=(u,\mathrm i u,v,\mathrm i v)$ with
$u,v$ independent and all other vectors zero; use
$P(\mathrm i u)=P(u)^{-1}$.  In either case the vectors span $V$.
The locus where a chosen pair has nonzero determinant has closed
free orbits and admits the usual frame quotient, of relative
dimension three.  These observations prove
$\dim\mathbb C[Z_{g,p}]^G\geq D-3$ in every case.

A nonzero ideal in the finitely generated domain $S_{g,p}$ lowers
dimension.  The surjection \eqref{RYevaluationmap} therefore has
zero kernel, proving the lemma.
\end{proof}

\begin{proof}[Proof of Theorem \ref{RYnormality}]
The invariant ring of a normal affine scheme under a reductive
group is normal.  Lemmas \ref{RYmatrixnormal} and
\ref{RYinvariantmodel} therefore prove that $S_{g,p}$ is normal.
Let
\[
 L=\mathbb C[v_1^{\pm1},\ldots,v_p^{\pm1}],\qquad
 L'=\mathbb C[t_1^{\pm1},\ldots,t_p^{\pm1}],\qquad v_j=t_j^2.
\]
The coefficient change in
\cite[Proposition 3.8]{bloomquist23degenerations} gives
\[
 B_{g,p}\otimes_L L'\simeq S_{g,p}\otimes_{\mathbb C}L'.
\]
The right side is normal.  Since $L'/L$ is finite, faithfully flat,
and \'etale, normality descends to $B_{g,p}$.  Finite generation and
the domain property follow from
\cite[Theorem 2]{bloomquist23degenerations}.
\end{proof}

\begin{prop}\label{RYspherefour}
The algebra $B_{0,4}$ is not normal.  In particular,
\[
 B_{0,4}\ \subsetneq\
 \Gamma(\mathring{\mathcal A},\mathcal O_{\mathring{\mathcal A}})
 =\overline{B_{0,4}},
\]
where the bar denotes integral closure in the common fraction field.
\end{prop}

\begin{proof}
Choose the tetrahedral ideal triangulation and write its six lambda
lengths as
\[
 \begin{gathered}
 a=x_{12},\quad b=x_{13},\quad c=x_{14},\quad
 d=x_{23},\quad e=x_{24},\quad f=x_{34},\\
 N=af+be+cd.
 \end{gathered}
\]
The four potential functions in this chart are
\[
 v_1=\frac N{abc},\qquad v_2=\frac N{ade},\qquad
 v_3=\frac N{bdf},\qquad v_4=\frac N{cef}.
\]
Consequently $h=a/f$ belongs to $\operatorname{Frac}(B_{0,4})$ and
satisfies the monic equation
\begin{equation}\label{RYspherefourintegral}
 h^2=\frac{v_3v_4}{v_1v_2}.
\end{equation}
To show that $h\notin B_{0,4}$, use the sphere presentation
\cite[Theorem 1.1 and formulas (4.6)--(4.7)]{azad21presentations}.
It admits the augmentation
\[
 \epsilon:B_{0,4}\longrightarrow
 L=\mathbb C[v_1^{\pm1},v_2^{\pm1},v_3^{\pm1},v_4^{\pm1}],
 \qquad \beta_{ij}\longmapsto0,\quad v_i\longmapsto v_i.
\]
Indeed, at $q^{1/2}=1$ the Ptolemy, commutation, and
$\gamma$-relations vanish under this substitution.  The constant
term of the big-circle polynomial is $2(-1)^{4-1}=-2$, as required
by the big-circle relation.  If $h$ belonged to $B_{0,4}$,
\eqref{RYspherefourintegral} would give a square root in $L$ of
$v_1^{-1}v_2^{-1}v_3v_4$.  This is impossible: an element whose
square is a Laurent unit is itself a Laurent unit, and the exponent
of every variable in its square is even.  Thus $B_{0,4}$ is not
integrally closed.  The last equality in the proposition is Proposition
\ref{normalizationcomparison}.  Alternatively, the smoothness of
$\mathring{\mathcal A}$ shows directly that the integral rational
function $h$ is regular there.
\end{proof}

\subsection{The Roger--Yang algebra and decorated Teichm\"uller space}

Roger and Yang constructed the lambda-function homomorphism
\[
 Sk^{RY}(\Sigma)\longrightarrow C^\infty(\mathcal T^d(\Sigma))
\]
and asked for an algebraic interpretation of decorated Teichm\"uller
space compatible with this algebra \cite{roger14skein}.  Injectivity
was established in \cite{moon24consequences,bloomquist23degenerations}.
We now distinguish the positive real locus from the entire real
point set, and use the independent normality theorem above.

The positive cluster transition maps define the real cluster scheme
$\mathcal A_{\mathbb R}$.  Let
$\mathring{\mathcal A}_{\mathbb R}$ be its largest open subscheme on
which every $v_i$ and $v_i^{-1}$ is regular.  Write
$\mathcal A_{\mathbb R}(\mathbb R_{>0})$ for the locus with positive
coordinates in the positive cluster atlas.  This is notation for a
set of positive real points; it does not denote morphisms from the
spectrum of $\mathbb R_{>0}$.

\begin{cor}\label{localized}
For a triangulable $\Sigma_{g,p}$ with $p\geq2$, the positive locus
satisfies
\[
 \mathring{\mathcal A}_{\mathbb R}(\mathbb R_{>0})
 =\mathcal A_{\mathbb R}(\mathbb R_{>0})
 \simeq\mathcal T^d(\Sigma_{g,p}).
\]
If, in addition, $g\geq1$ or $p\geq5$, then
\[
 Sk^{RY}(\Sigma_{g,p})
 =\Gamma(\mathring{\mathcal A},\mathcal O_{\mathring{\mathcal A}}),
 \qquad
 Sk^{RY}(\Sigma_{g,p})_{\mathbb R}
 =\Gamma(\mathring{\mathcal A}_{\mathbb R},
          \mathcal O_{\mathring{\mathcal A}_{\mathbb R}}).
\]
For $\Sigma_{0,4}$, the complex global-section ring is instead the
strict normalization described in Proposition \ref{RYspherefour}.
\end{cor}

\begin{proof}
The identification of the positive cluster locus with decorated
Teichm\"uller space is
\cite[Theorem 1.7(b)]{fock06moduli}, with the usual compatible lift
convention.  The potentials are positive on this locus, so in every
positive chart each $v_i$ and $v_i^{-1}$ is regular at each such
point.  This proves the assertion about positive points.

The complex equality follows from Proposition
\ref{normalizationcomparison} and Theorem \ref{RYnormality}.
For the real equality, the construction and the common open
immersion of Proposition \ref{commonlocus} are defined over
$\mathbb R$.  In particular,
$\mathring{\mathcal A}_{\mathbb R}$ is a quasi-compact separated
open subscheme of the affine scheme
$\operatorname{Spec}(Sk^{RY}(\Sigma)_{\mathbb R})$.
Using a finite affine cover and the equalizer computing global
sections, flat base change gives
\[
 \Gamma(\mathring{\mathcal A}_{\mathbb R},\mathcal O)
       \otimes_{\mathbb R}\mathbb C
 \simeq\Gamma(\mathring{\mathcal A},\mathcal O).
\]
The real skein algebra also complexifies to the complex skein
algebra.  Therefore the real comparison becomes an isomorphism
after tensoring with $\mathbb C$ and is already an isomorphism by
faithful flatness.
\end{proof}

\subsection{Algebraic properties of skein algebras}

\begin{prop}
For a triangulable $\Sigma_{g,p}$ with $p\geq2$, the tagged skein
algebra is a finitely generated normal Gorenstein domain.
The Roger--Yang algebra is a finitely generated domain, and is
normal when $g\geq1$ or $p\geq5$; it is not normal for $\Sigma_{0,4}$.
\end{prop}

\begin{proof}
Finite generation and the Gorenstein property of the tagged skein
algebra are Theorem \ref{skeinalgebraic}.  By Theorem
\ref{mainthmproof}, it equals $U(\Sigma)$, an intersection of Laurent
polynomial rings in their common fraction field, and hence is
normal.  For the Roger--Yang algebra, finite generation and the
domain property are \cite[Theorem 2]{bloomquist23degenerations};
normality and the exception follow from Theorem \ref{RYnormality}
and Proposition \ref{RYspherefour}.
\end{proof}

\begin{rem}
The normality proof for the Roger--Yang algebra is independent of
the Gorenstein property of the tagged skein algebra.  No
Cohen--Macaulay or Gorenstein assertion for the Roger--Yang algebra
is needed for the global-section comparison.
\end{rem}

\subsection{Finite generation of \texorpdfstring{$U(\Sigma)$}{the upper cluster algebra}}

\begin{prop}
For a triangulable surface $\Sigma$ covered by the hypotheses of
Corollary \ref{borderedAskeinU}, or for a triangulable surface
$\Sigma_{g,p}$ without boundary, the upper cluster algebra
$U(\Sigma)$ is finitely generated over $\mathbb C$.
\end{prop}

\begin{proof}
In the boundary case, use Corollary \ref{borderedAskeinU} and finite
generation of the tagged skein algebra.  For $p\geq2$, use Theorems
\ref{mainthmproof} and \ref{skeinalgebraic}.  For $p=1$, Proposition
\ref{onceuppernormalization} identifies $U(\Sigma_{g,1})$ with the
normalization of the finitely generated domain
$Sk^{ta}(\Sigma_{g,1})[v]$.  The normalization of a finitely
generated domain over $\mathbb C$ is finite over that domain, so
$U(\Sigma_{g,1})$ is finitely generated as well.
\end{proof}

\section{Once-punctured surfaces without boundary}
For $\Sigma=\Sigma_{g,1}$, $g\geq 1$, the ordinary cluster variables correspond to the plain arcs, or, after changing the tag at the puncture, to the notched arcs (Theorem \ref{arcvariable}). These two collections belong to different components of the tagged exchange graph. We therefore compare the tagged skein algebra with the union of the two cluster atlases.

Write $R=Sk^{ta}(\Sigma)$ and let $\mathcal A$ be the cluster variety associated with the plain arcs. Changing the tag at the puncture gives a second cluster variety $\mathcal A'$ and a birational map $\mathcal A\dashrightarrow\mathcal A'$. Glue the two varieties along the largest open subsets on which this map is an isomorphism, and denote the resulting variety by $\widetilde{\mathcal A}$. Every tagged arc and loop is Laurent in every chart of either atlas. Consequently, if $x_1,\ldots,x_d$ are the variables of such a chart, then
\[
 R[(x_1\cdots x_d)^{-1}]
 =\mathbb C[x_1^{\pm1},\ldots,x_d^{\pm1}],
 \qquad d=6g-3.
\]
Thus the cluster charts, glued along their maximal overlaps, define an open immersion
$\widetilde{\mathcal A}\hookrightarrow\operatorname{Spec}R$.

\begin{thm}
For $\Sigma=\Sigma_{g,1}$, $g\geq 1$, the complement of
$\widetilde{\mathcal A}\hookrightarrow\operatorname{Spec}(Sk^{ta}(\Sigma))$
has codimension at least two. In particular, $Sk^{ta}(\Sigma)$ satisfies R1, and
\[
 \Gamma(\widetilde{\mathcal A},\mathcal O_{\widetilde{\mathcal A}})
 =\overline{Sk^{ta}(\Sigma)},
\]
where the bar denotes integral closure in the common fraction field. If $Sk^{ta}(\Sigma)$ also satisfies S2, then
$Sk^{ta}(\Sigma)=\Gamma(\widetilde{\mathcal A},\mathcal O_{\widetilde{\mathcal A}})$.
\label{oncepuncturedconj}
\end{thm}

\begin{proof}
Let $I\subset R$ be the ideal generated by all tagged arcs. By the puncture-ideal estimate in Lemma \ref{dimensionlemma},
\[
 \dim(R/I)\leq d-3,
 \qquad d=\dim R=6g-3.
\]
In particular, every prime of height one in $R$ lies in $D(a)$ for some tagged arc $a$.

Choose the component of the cluster atlas containing $a$, and freeze $a$ in a seed containing it. Let $\mathcal V_a$ be the union of the cluster charts obtained without mutating $a$, and let $A^{[a]}$ and $U^{[a]}$ be the corresponding frozen cluster and upper cluster algebras, with $a$ inverted. Universal Laurentness gives
\[
 A^{[a]}\ \subseteq\ R[a^{-1}]\ \subseteq\ U^{[a]}.
\]
By Lemma \ref{onearcfreezing}, the two outer algebras are equal, including when cutting along $a$ produces two unpunctured components. Hence
\[
 R[a^{-1}]=U^{[a]}
 =\Gamma(\mathcal V_a,\mathcal O_{\mathcal V_a}).
\]
The cluster charts define an open immersion
$\mathcal V_a\hookrightarrow\operatorname{Spec}R[a^{-1}]$, and Lemma \ref{basiclemma} shows that its complement has codimension at least two. Since $\mathcal V_a\subseteq\widetilde{\mathcal A}$, every height-one point of $D(a)$ belongs to $\widetilde{\mathcal A}$. The estimate for $I$ therefore proves the asserted codimension bound. As $\widetilde{\mathcal A}$ is smooth, it also proves R1 for $R$.

The ring $R$ is a finitely generated domain over $\mathbb C$, so its normalization $\overline R$ is finite over $R$. The normalization morphism is an isomorphism over the normal open subset $\widetilde{\mathcal A}$. Its inverse image has complement of codimension at least two: for finite morphisms of affine varieties, dimensions of closed subsets are preserved. The S2 property of $\overline R$ now gives
\[
 \overline R=\Gamma(\widetilde{\mathcal A},\mathcal O_{\widetilde{\mathcal A}}).
\]
Finally, if $R$ satisfies S2, then R1 and Serre's criterion imply $R=\overline R$.
\end{proof}

The usual upper cluster algebra has a similar description after adjoining the potential. This description also supplies the once-punctured case of finite generation in Section 6.

\begin{prop}
For $\Sigma=\Sigma_{g,1}$, $g\geq 1$, let $v$ be the potential regular on the plain cluster atlas $\mathcal A$. Then
\[
 U(\Sigma)=\overline{Sk^{ta}(\Sigma)[v]},
\]
where integral closure is taken in the cluster function field. In particular, $U(\Sigma)$ is finitely generated over $\mathbb C$.
\label{onceuppernormalization}
\end{prop}

\begin{proof}
Put $C=R[v]\subseteq U(\Sigma)$, and let $J\subset C$ be the ideal generated by all plain arcs. Both $R$ and $C$ have fraction field equal to the cluster function field and dimension $d=6g-3$. If $a$ is a plain arc, its all-notched version is $v^2a$. Thus $IC\subseteq J$, where $I$ is the ideal from the preceding proof. The natural map $R[T]\to C$, $T\mapsto v$, induces a surjection
\[
 (R/I)[T]\longrightarrow C/J.
\]
Consequently,
\[
 \dim(C/J)\leq\dim(R/I)+1\leq d-2.
\]

For every plain arc $a$, the generalized tagged arc with one notched end represents $va$ and belongs to $R$ by Definition--Proposition \ref{taggedskein}. Therefore $v\in R[a^{-1}]$ and
\[
 C[a^{-1}]=R[a^{-1}]=U^{[a]}.
\]
As in the preceding proof, the frozen cluster charts in the plain component cover all height-one points of $\operatorname{Spec}C[a^{-1}]$. Moreover, each ordinary cluster chart is a principal open subset of $\operatorname{Spec}C$, since $C\subseteq U(\Sigma)$ contains its cluster variables. We obtain an open immersion
$\mathcal A\hookrightarrow\operatorname{Spec}C$ whose complement has codimension at least two, using the dimension estimate for $C/J$.

Normalize $C$. The normalization is finite because $C$ is a finitely generated domain over $\mathbb C$. It is an isomorphism over the smooth open subset $\mathcal A$, whose complement still has codimension at least two. Taking global sections and using normality gives
\[
 \overline C=\Gamma(\mathcal A,\mathcal O_{\mathcal A})=U(\Sigma).
\]
The finiteness of the normalization proves the final assertion.
\end{proof}

\begin{rem}
An enlarged cluster variety for the once-punctured torus is also studied using scattering diagrams in \cite{zhou20cluster}. Identifying that construction with the maximal-overlap gluing used here requires a separate comparison. The theorem above does not rely on that identification. The equality
$Sk^{ta}(\Sigma)=\Gamma(\widetilde{\mathcal A},\mathcal O_{\widetilde{\mathcal A}})$
for all $g\geq1$ remains conditional on S2 for $Sk^{ta}(\Sigma)$.
\end{rem}

\end{document}